\documentclass[10pt]{amsart}

\usepackage{IEEEtrantools}
\usepackage{amsfonts}
\usepackage{amssymb}
\usepackage{amsmath,amsfonts,amsthm}
\usepackage{mathrsfs}
\usepackage{url}
\usepackage{color,soul}

\newcommand{\scrt}{\mathscr{T}}

\newcommand{\trace}{\mbox{\textit{Trace}\,}}

\newcommand{\be}{\begin{equation}}
\newcommand{\ee}{\end{equation}}
\newcommand{\bea}{\begin{eqnarray}}
\newcommand{\eea}{\end{eqnarray}}
\newcommand{\bean}{\begin{eqnarray*}}
\newcommand{\eean}{\end{eqnarray*}}
\newcommand{\brray}{\begin{array}}
\newcommand{\erray}{\end{array}}
\newcommand{\biearray}{\begin{IEEEarray}{rCl}}
\newcommand{\eiearray}{\end{IEEEarray}}

\newtheorem{dfn}{Definition}[section]
\newtheorem{thm}[dfn]{Theorem}
\newtheorem{lmma}[dfn]{Lemma}
\newtheorem{ppsn}[dfn]{Proposition}
\newtheorem{crlre}[dfn]{Corollary}
\newtheorem{xmpl}[dfn]{Example}
\newtheorem{rmrk}[dfn]{Remark}

\newcommand{\bdfn}{\begin{dfn}\rm}
\newcommand{\bthm}{\begin{thm}}
\newcommand{\blmma}{\begin{lmma}}
\newcommand{\bppsn}{\begin{ppsn}}
\newcommand{\bcrlre}{\begin{crlre}}
\newcommand{\bxmpl}{\begin{xmpl}}
\newcommand{\brmrk}{\begin{rmrk}\rm}

\newcommand{\edfn}{\end{dfn}}
\newcommand{\ethm}{\end{thm}}
\newcommand{\elmma}{\end{lmma}}
\newcommand{\eppsn}{\end{ppsn}}
\newcommand{\ecrlre}{\end{crlre}}
\newcommand{\exmpl}{\end{xmpl}}
\newcommand{\ermrk}{\end{rmrk}}

\newcommand{\bbc}{\mathbb{C}}
\newcommand{\bbz}{\mathbb{Z}}
\newcommand{\bbn}{\mathbb{N}}

\newcommand{\cla}{\mathcal{A}}
\newcommand{\clb}{\mathcal{B}}
\newcommand{\clh}{\mathcal{H}}

\newcommand{\cll}{\mathcal{L}}

\newcommand{\cls}{\mathcal{S}}

\newcommand{\prf}{\noindent{\it Proof\/}: }

\def \qed { \mbox{}\hfill
$\Box$\vspace{1ex}}

\begin{document}
\title[Local index formula for the quantum double suspension]{Local index formula for the quantum double suspension}
\author{Partha Sarathi Chakraborty}
\address{The Institute of Mathematical Sciences, CIT Campus, Taramani, Chennai
600113}
\email{parthac@imsc.res.in}
%\thanks{Partha Sarathi Chakraborty acknowledges financial support from Indian National Science Academy through its project ``Noncommutative Geometry of Quantum Groups''}
\author{Bipul Saurabh}
\address{Indian Statistical Institute, 7 S.J.S.S Marg, NewDelhi 100016}
\email{bipul9r@isid.ac.in}
\thanks{First author  acknowledges support from Swarnajayanthi Fellowship  Award Project No. DST/SJF/MSA-01/2012-13\\ Corresponding author: Bipul Saurabh}
%\thanks{Corresponding author}

%
%\title{Local index formula for the quantum double suspension}
%
%\author{Partha Sarathi Chakraborty \corref{cor1} \fnref{fn1}}
%%\ead{parthac@imsc.res.in}
%\address{The Institute of Mathematical Sciences, CIT Campus, Taramani, Chennai
%600113}
%\email{parthac@imsc.res.in}
%\thanks{Partha Sarathi Chakraborty acknowledges financial support from Swarnajayanthi Fellowship  Award Project No. DST/SJF/MSA-01/2012-13}
%\author{Bipul Saurabh}
%\email{saurabhbipul2@gmail.com}
%\address{The Institute of Mathematical Sciences, CIT Campus, Taramani, Chennai
%600113}
%%\cortext[cor1]{Corresponding author}
%%\fntext[fn1]{First author  acknowledges support from Swarnajayanthi Fellowship  Award Project No. DST/SJF/MSA-01/2012-13\\ Corresponding author: Partha Sarathi Chakraborty}
%%\thanks{Corresponding author}
%\keywords{Spectral triples,  noncommutative geometry, Local Index Formula, Quantum Double Suspension,  Noncommutative Torus}
%%\MSC Primary 46L87 \sep 58B34 \sep  Secondary 19K56\end{keyword}
%\subjclass[2000]{Primary 46L87,58B34, Secondary 19K56}
\begin{abstract} 
Our understanding of local index formula in noncommutative geometry is stalled for a while because we do not
have more than one explicit computation, namely that of Connes for quantum SU(2) and do not understand the meaning of the various multilinear functionals involved in the formula. 
In such a situation further progress in understanding necessitates more explicit computations and here we  execute the second explicit computation for the quantum double suspension, a 
construction inspired by the Toeplitz extension. More specifically we compute local index formula for the quantum double suspensions of $C(\cls^2)$ and the noncommutative $2$-torus.

%In all these years we have seen essentially one explicit computation of the local index formula, namely that of quantum SU(2) done by Connes. Here we execute the second explicit computation. More specifically we compute 
% local index formula for the quantum double suspensions of $C(\cls^2)$ and the noncommutative $2$-torus.
%
	
\end{abstract}

\date{\today}

%%%%%%%%%%%%%%%%%%%%%%%%%%%%%%%%%

%%%%%%%%%%%%%%%%%%%%%%%%%%%%%%%%%%
%%%%%  ABSTRACT
%%%%%%%%%%%%%%%%%%%%%%%%%%%%%%%%%%
 
 \maketitle

% {\bf AMS Subject Classification No.:} {\large 58}B{\large 34}, {\large
% 46}L{\large 87}, {\large
%   19}K{\large 33}\\
% {\bf Keywords.} Spectral triples, noncommutative geometry, Local Index Formula, Quantum Double Suspension, Noncommutative Torus.
\section{INTRODUCTION}
In the spectral formulation of noncommutative geometry, Connes (\cite{Con-1994aa}) specifies a noncommutative geometric `space' by a 
triple consisting of  a Hilbert space ${\mathcal H}$, an involutive subalgebra $\mathcal A$ of the algebra of bounded operators on 
$\mathcal H$ and a self adjoint operator $D$ with compact resolvent. The algebra $\mathcal A$ and the operator 
$D$ are tied up by the requirement that the commutators $[D,{\mathcal A}]$ give rise to bounded operators. 
Such a triple is called a spectral triple or  an unbounded K-cycle. Often it is required that the spectral triple satisfies further 
 conditions (\cite{Con-2013aa}). The conditions of regularity and discrete dimension spectrum was introduced by 
 Connes and Moscovici in \cite{ConMos-1995aa}.
 More specifically a spectral triple $({\mathcal A}, {\mathcal H},D)$ is said to be regular if both $\mathcal A$ and $[D,{\mathcal A}]$ are 
 in the domains of  $\delta^n$ for all $n \ge 0$, where $\delta$ is the derivation $[|D|,\cdot]$. One says that the spectral triple 
 has  dimension spectrum ${\mathfrak S}$, if for every element $b$ in the smallest algebra $\mathcal B$ containing
 $\mathcal A$, $[D,{\mathcal A}]$ and closed under the derivation $\delta$, the associated zeta function $\zeta_b(z)=Tr b {|D|}^{-z}$ a priori 
 defined on the right half plane $\Re (z) > p$ admits a meromorphic extension to the  whole complex plane
 with poles contained in $\mathfrak S$. They arrived at these conditions in their efforts to give an expression for the Chern character in 
 terms of 
local data. The local Chern character constructed in (\cite{ConMos-1995aa})
is a cocycle in the total complex of the b-B bicomplex and therefore 
given by a sum of several multilinear functionals $\phi_n$'s. At this point one sees a departure from 
the classical world in the sense that for the canonical spectral triple associated with a Riemannian spin 
manifold most of these terms vanish( (b) page 231 in (\cite{ConMos-1995aa})). 
This indicates that the terms $\{\phi_n\}$ appearing in the local Chern character are offshoots of noncommutativity and may exhibit 
features of the noncommutative world not visible in the commutative side of the story. This calls for further exploration of 
the local Chern character. 

Indeed during 2000-2006 we saw a lot of activity around the local index formula (LIF). Most of these (\cite{CarPhiRenSuk-2006aa},\cite{CarPhiRenSuk-2006ab},\cite{CarPhiRenSuk-2013aa}) were concerned with 
extending the  formula itself to the framework of semi-finite spectral triples, while Connes (\cite{Con-2004aa}) gave a demonstration of the formula for the spectral triple constructed in (\cite{PsPal_SU_q_2}).
This computation of Connes is our point of departure.
One might wonder what is the point of the article (\cite{Con-2004aa})?  Should one interpret it as just a one time calculation or as an invitation to make explicit computations in various instances to gain
further  insight on the  terms appearing in the formula. %This line of exploration is also in line with  the 
%general principle outlined in \cite{Con-Interview}. 
We take the second interpretation. But unfortunately in all these years we have only one  (\cite{VanDab-2005aa}) more explicit computation of the local index formula. 
That too essentially follows from the arguments of Connes  and  the terms that contribute to the residue cocycle  in (\cite{VanDab-2005aa}) are not different from those found by Connes.
Therefore, so far, effectively we have only one explicit computation of the local index formula. In another attempt, (\cite{PalSun-2010aa}) did a careful analysis of (\cite{Con-2004aa})
and extended the proof of regularity and discreteness of dimension spectrum  to the case of odd dimensional quantum spheres, but did not quite compute the formula. In such a situation it
looked imperative to try and compute  in some more cases. The first stumbling block is establishing regularity and discreteness of the dimension spectrum.
%
%
%
%
%To understand the real content of these functionals one might attempt to compute these in 
%examples. This is really step one of a two step strategy described in (\cite{Con-Interview}). 
%In that endeavor one soon realizes that in the examples studied by Connes and Moscovici (\cite{ConMos-1995aa}) there are too 
%many terms and they cancel.  This was explained by them in (\cite{ConMos-1998aa}). The first `tractable' instance of a 
%computation was demonstrated by Connes in (\cite{Con-2004aa}). 
%However for quite sometime there were no further progress. 
%The situation changed 
That was achieved in (\cite{ChaSun-2011ab}), through the   construction of 
quantum double suspension (QDS) of spectral triples hitherto known for $C^*$-algebras (\cite{HonSzy-2002aa}). 
Chakraborty and Sundar also identified the hypothesis of weak heat kernel 
asymptotic expansion (to be abbreviated  as WHKAE) that allows one to conclude stability of the hypothesis of regularity and dimension spectrum under QDS. 
More specifically it was shown that if one 
starts with a spectral triple that is regular and has got discrete dimension spectrum satisfying WHKAE then its quantum double 
suspension is regular, has discrete dimension spectrum and 
also satisfies WHKAE. The computation in (\cite{Con-2004aa}) can also be seen in this light as the QDS of the canonical spectral 
triple of the circle. But (\cite{ChaSun-2011ab}) fell short of the actual
description of the multilinear functionals involved in the local Chern character and that 
brings us to the contents of this paper, namely, explicit computation of the LIF for the quantum double suspension.  
In the terminology of (\cite{ChaSun-2011ab}) the articles
(\cite{Con-2004aa},\cite{PalSun-2010aa}) considered LIF  for $\Sigma^2 C(S^1)$, the quantum double suspension  of the 
circle and its iterations. Here we consider LIF for the quantum double suspension of 
two dimensional manifolds. More specifically we take up two cases, one classical and one noncommutative. The classical case we consider is 
that of the two sphere and the noncommutative
case tackled is the noncommutative two torus. Thus this article can be seen  as the second computation  after that of (\cite{Con-2004aa}). 
In section two we recall the local Chern character, the main object of interest in the LIF. 
It is a cocycle in the total complex of the $b-B$-bicomplex and depending upon the parity of 
the spectral triple is given by a finite sequence of multilinear functionals $\{\phi_{2n}\}$ or $\{\phi_{2n+1}\}$. 
For our purpose we express $\phi_n$'s in terms of some other functional
$\psi_x^{(0)}$'s. Let $\Sigma^2\psi_x^{(0)}$ denote the corresponding quantity for the quantum double suspended spectral triple. In 
section three we obtain expressions for these. Finally
in section four we apply these in two concrete situations. We take up the case of quantum double suspension of the two sphere
and the QDS of the noncommutative two torus. 
We explicitly describe $\Sigma^2\phi_n$'s, the multilinear functionals involved in the LIF. Section five contains concluding remarks.

\section{LOCAL INDEX FORMULA}
Let $(\cla, \clh ,D)$ be a $p^+$ summable regular spectral triple. We denote by $d$ , $\delta$ and $\nabla$ the commutator
with $D$, $|D|$ and $D^2$ respectively 
i.e. $da=\left[D,a\right], \delta a =\left[\left|D\right|,a\right]$, and $\nabla a=\left[D^2, a\right]$.
For $k \in \bbn^r$, define $\left|k\right| =: \sum_{i=1}^rk_i$.  
We will now state the Connes--Moscovici local index theorem for a regular odd spectral triple with discrete dimension spectrum.

%%%%%%%%%%%%%%%%%%%%%%%%%%%%%%%%%%%%%%%%%%%
\bthm [\cite{ConMos-1995aa}]
Let $(\cla,\clh,D)$ be a regular $p^+$ summable odd spectral triple.
Assume further that the dimension spectrum is discrete with
finite multiplicity.
Define, for $n$ odd,
\begin{eqnarray*}
&& \phi_n(a_0,\ldots,a_n)\\
&=&\sum_{k_j\geq 0}
  c_{n,k} \mbox{Res}_{z=0}\left(\Gamma(|k| +\frac{n}{2}+z) \trace \left(a_0\nabla^{k_1}(da_1)\ldots \nabla^{k_n}(da_n)|D|^{-n-2|k|-2z}\right)\right),
\end{eqnarray*}
where $a_j\in\cla$, $k=(k_1,\ldots,k_n)\in\bbn^n$ and $c_{n,k}$ are given by
\[
 c_{n,k}=\sqrt{2i}(-1)^{|k|} (\prod_{j=1}^n k_j! \prod_{j=1}^n  (k_1+\ldots +k_j+j))^{-1}
\]
Then $\phi_n$ is zero except for finitely many $n$'s and 
$(\phi_1,\phi_3,\ldots)$ is a $(b,B)$-cocycle.

Furthermore, the cohomology class of this cocycle in $HC^{odd}(\cla)$ is same as the
Chern character of $(\cla,\clh,D)$,  in particular, for $[u]\in K_1(\cla)$, one has
\[
 \langle(\phi_n),[u]\rangle=\langle [Ch^F],[u]\rangle.
\]
\ethm 
%%%%%%%%%%%%%%%%%%%%%%%%%%%%%%%%%%%%%%%%%%%
%For a proof, see the paper  (\cite{con-mos-1995a}) by Connes and Moscovici.

%%%%%%%%%%%%%%%%%%%%%%%%%%%%%%%%%%%%%%%%%%%
\bcrlre\label{cor:localindex}
If in addition we assume that the dimension spectrum is discrete and simple,
then the cocycle $\phi_n$ in the above theorem is given by
\[
 \phi_n(a_0,\ldots,a_n):=\sum_k c_{n,k}\mbox{Res}_{z=0} 
   \trace\left(a_0\nabla^{k_1}(da_1)\ldots \nabla^{k_n}(da_n)|D|^{-n-2|k|-2z}\right),
\]
where $a_j\in\cla$, $k=(k_1,\ldots,k_n)\in\bbn^n$ and $c_{n,k}$ are given by
\[
 c_{n,k}=(-1)^{|k|}\sqrt{2i} (\prod k_j! \prod (k_1+\ldots +k_j+j))^{-1}
   \Gamma(|k| +\frac{n}{2}).
\]
\ecrlre 
For a regular even spectral triple with discrete dimension spectrum, we state the 
Connes--Moscovici local index theorem as follows :
\bthm[\cite{Hig-2004aa}]
Let $(\cla,\clh,D,\gamma)$ be a regular $p^+$ summable even spectral triple.
Assume further that the dimension spectrum is discrete with
finite multiplicity.
Define, for $n$ even,
\begin{IEEEeqnarray}{rCl}
\phi_0(a_0) &=& \mbox{Res}_{z=0}\left(\Gamma(z)\trace(\gamma a_0\left|D\right|^{-2z})\right). \nonumber \\
\phi_n(a_0,\cdots,a_n)&=&\sum_{k_j\geq 0} c_{n,k} \mbox{Res}_{z=0}\left(\Gamma(|k| +\frac{n}{2}+z) \trace \left(\gamma a_0\nabla^{k_1}(da_1)\ldots \nabla^{k_n}(da_n)|D|^{-n-2|k|-2z}\right)\right). \nonumber
\end{IEEEeqnarray}
where $a_j\in\cla$, $k=(k_1,\ldots,k_n)\in\bbn^n$ and $c_{n,k}$ are given by
\[
 c_{n,k}=\sqrt{2i}(-1)^{|k|} (\prod_{j=1}^n k_j! \prod_{j=1}^n  (k_1+\ldots +k_j+j))^{-1}
\]
Then $\phi_n$ is zero except for finitely many $n$'s and 
$(\phi_0,\phi_2,\ldots)$ is a $(b,B)$-cocycle.

Furthermore, the cohomology class of this cocycle in $HC^{even}(\cla)$ is same as the
Chern character of $(\cla,\clh,D,\gamma)$, in particular,
for $[p]\in K_0(\cla)$, one has
\[
 \langle(\phi_n),[p]\rangle=\langle [Ch^F],[p]\rangle.
\]
\ethm 
%%%%%%%%%%%%%%%%%%%%%%%%%%%%%%%%%%%%%%%%%%%
%For a proof, see the paper  (\cite{con-mos-1995a}) by Connes and Moscovici.

%%%%%%%%%%%%%%%%%%%%%%%%%%%%%%%%%%%%%%%%%%%
\bcrlre\label{cor: evenlocalindex }
If in addition we assume that the dimension spectrum is discrete and simple,
then the cocycle $\phi_n$ in the above theorem is given by
\begin{IEEEeqnarray}{rCl}
\phi_0(a_0)&=& \mbox{Res}_{z=0}\left(\Gamma(z)\trace\left(\gamma a_0\left|D\right|^{-2z}\right)\right). \nonumber \\
\phi_n(a_0,\ldots,a_n)&=&\sum_k c_{n,k}\mbox{Res}_{z=0}\trace\left(a_0\nabla^{k_1}(da_1)\ldots \nabla^{k_n}(da_n)|D|^{-n-2|k|-2z}\right), \nonumber
\end{IEEEeqnarray}
where $a_j\in\cla$, $k=(k_1,\ldots,k_n)\in\bbn^n$ and $c_{n,k}$ are given by
\[
 c_{n,k}=(-1)^{|k|}\sqrt{2i} (\prod k_j! \prod (k_1+\ldots +k_j+j))^{-1}
   \Gamma(|k| +\frac{n}{2}).
\]
\ecrlre 
By rearranging the terms, we  write the linear functionals $\phi_n$ as sum of some
other linear functionals which will be more tractable for our purpose.   Note that
\begin{IEEEeqnarray}{rCl}
\nabla^n(T)& =&\sum_{k=0}^n 2^{n-k}{n\choose k}\delta^{n+k}(T)|D|^{n-k}. \nonumber \\
 \left|D\right|^nT& =&\sum_{k=0}^n {n\choose k}\delta^{k}(T)|D|^{n-k}. \nonumber 
\end{IEEEeqnarray} 
Using these equations, we get
\begin{IEEEeqnarray}{lCl}
 \nabla^{k_1}(T_1)\nabla^{k_2}(T_2) \nonumber \\
 =\sum_{j_1=0}^{k_1} \sum_{j_2=0}^{k_2}2^{k_1-j_1}2^{k_2-j_2}{k_1\choose j_1}{k_2\choose j_2}\delta^{k_1+j_1}(T_1)|D|^{k_1-j_1}
 \delta^{k_2+j_2}(T_2)|D|^{k_2-j_2}. \nonumber \\
=\sum_{j_1=0}^{k_1} \sum_{j_2=0}^{k_2} \sum_{s_1=0}^{k_1-j_1}2^{k_1-j_1}2^{k_2-j_2}{k_1\choose j_1}{k_2\choose j_2}{k_1-j_1\choose s_1}
\delta^{k_1+j_1}(T_1)\delta^{k_2+j_2+s_1}(T_2)|D|^{k_1-j_1+k_2-j_2-s_1}. \nonumber 
\end{IEEEeqnarray}
Hence
\begin{IEEEeqnarray}{rCl}
 a_{0} \nabla^{k_1}(a_1)\nabla^{k_2}(a_2) \cdots \nabla^{k_n}(a_n) \left|D\right|^{-n-2\left|k\right|-2z}&=&\sum_{j_i=0}^{k_i} \sum_{s_1=0}^{S_1^{'}}\sum_{s_2=0}^{S_2^{'}}\cdots \sum_{s_n=0}^{S_n^{'}}2^{\left|k\right|-\left|j\right|}\prod_{i=1}^n{k_i\choose j_i}{S_i^{'}\choose s_i}a_0 \delta^{k_1+j_1}(a_1)\nonumber \\
& & \delta^{k_2+j_2+s_1}(a_2)\cdots \delta^{k_n+j_n+s_{n-1}}(a_n)\left|D\right|^{-n-\left|k\right|-\left|j\right|-\left|s\right|-2z}. \nonumber 
\end{IEEEeqnarray}  
where $S_i^{'}=\sum_{\ell=1}^{i-1}(k_{\ell}-j_{\ell}-s_{\ell}) + k_i-j_i$. Define $x_1= k_1+j_1$ and for $i>1$, 
$x_i=k_i+j_i+s_{i-1}$. 
Reparametrizing in the variables $x,k$ and $s$, we have
\begin{IEEEeqnarray}{lCl}
\sum_k c_{n,k}\left(a_0\nabla^{k_1}(da_1)\ldots \nabla^{k_n}(da_n)\left|D\right|^{-n-2\left|k\right|-2z}\right) \nonumber \\
= \sum_{x_i=0}^{\infty}\sum_{k_1=\left\lceil x_1/2\right\rceil}^{x_1}\sum_{s_1=0}^{S_1}\sum_{k_2=\left\lceil \left(x_2-s_1\right)/2\right\rceil}^{x_2}\sum_{s_2=0}^{S_2} \cdots \sum_{s_{n-1}=0}^{S_{n-1}}\sum_{k_n=\left\lceil \left(x_n-s_{n-1}\right)/2\right\rceil}^{x_n} \nonumber \\
c_{n,k} 2^{2\left|k\right|-\left|x\right|+\left|s\right| }\prod_{i=1}^n{k_i\choose x_i-k_i-s_{i-1}}{S_i\choose s_i}a_0\delta^{x_1}(a_1) \cdots \delta^{x_n}(a_n)\left|D\right|^{-n-\left|x\right|-2z}. \nonumber
 \end{IEEEeqnarray} 
where $S_i= \sum_{\ell=1}^{i} (2k_{\ell}-x_{\ell})$. Therefore
\begin{IEEEeqnarray}{lCl}
\phi(a_0,a_1, \cdots ,a_n)=\sum_{x_i=0}^{\infty}\sum_{k_1=\left\lceil x_1/2\right\rceil}^{x_1}\sum_{s_1=0}^{S_1}\sum_{k_2=\left\lceil \left(x_2-s_1\right)/2\right\rceil}^{x_2}\sum_{s_2=0}^{S_2} \cdots \sum_{s_{n-1}=0}^{S_{n-1}}\sum_{k_n=\left\lceil \left(x_n-s_{n-1}\right)/2\right\rceil}^{x_n} \nonumber \\
c_{n,k} 2^{2\left|k\right|-\left|x\right|+\left|s\right| }\prod_{i=1}^n{k_i\choose x_i-k_i-s_{i-1}}{S_i\choose s_i}\mbox{Res}_{z=0}\trace \left(a_0\delta^{x_1}(da_1) \cdots \delta^{x_n}(da_n)\left|D\right|^{-n-\left|x\right|-2z}\right). \nonumber \\
= \sum_{x_i=0}^{\infty}B_{x}^n\mbox{Res}_{z=\left(n+\left|x\right|\right)/2}\trace \left(a_0\delta^{x_1}(da_1) \cdots \delta^{x_n}(da_n)\left|D\right|^{-2z}\right). \nonumber
\end{IEEEeqnarray}  
where $B_x^n$ is given by
\begin{IEEEeqnarray}{lCl} 
B_x^n=\sum_{k_1=\left\lceil x_1/2\right\rceil}^{x_1}\sum_{s_1=0}^{S_1}\sum_{k_2=\left\lceil \left(x_2-s_1\right)/2\right\rceil}^{x_2}\sum_{s_2=0}^{S_2} \cdots \sum_{k_n=\left\lceil \left(x_n-s_{n-1}\right)/2\right\rceil}^{x_n}
c_{n,k} 2^{2\left|k\right|-\left|x\right|+\left|s\right| }\prod_{i=1}^n{k_i\choose x_i-k_i-s_{i-1}}{S_i\choose s_i} \nonumber 
\end{IEEEeqnarray}
Define for $x \in \bbn^n$ and $k \in \bbn$,
\begin{IEEEeqnarray}{rCl}
\psi_{x}^k(a_0,a_1, \cdots ,a_n)&=&\mbox{Res}_{z=\left(n+\left|x\right|+k\right)/2}\trace \left(a_0\delta^{x_1}(a_1) \cdots \delta^{x_n}(a_n)\left|D\right|^{-2z}\right) \nonumber 
\end{IEEEeqnarray}
Then $\phi_n$ can be written as
\begin{IEEEeqnarray}{rCl}
\phi_n(a_0,a_1, \cdots ,a_n)&=&\sum_{x_i=0}^{\infty}B_x^n\psi_x^{(0)}(a_0,da_1, \cdots ,da_n) \label{e4}
\end{IEEEeqnarray} 
Note that if $n+\left|x\right|+k>p$ then $\psi_{x}^k(a_0,a_1, \cdots ,a_n) =0$ for all  $(a_0,a_1, \cdots ,a_n)$. 
So, the sum given above is a finite sum. 
For an even spectral triple $(\cla, \clh,D,\gamma)$, one has 
\begin{IEEEeqnarray}{rCl}
\phi_n(a_0,a_1, \cdots ,a_n)&=&\sum_{x_i=0}^{\infty}B_x^n\psi_x^{(0)}(\gamma a_0,da_1, \cdots ,da_n) \label{e5}
\end{IEEEeqnarray}

\section{LOCAL INDEX FORMULA FOR QUANTUM DOUBLE SUSPENSION}

Let us fix some notations. We denote the left shift on $\ell^{2}\left(\bbn\right)$ by $S$ which is defined on the
standard orthonormal basis $\left(e_n\right)$ as $Se_n=e_{n-1}$. For $n < 0$, we denote by $S^n$ the operator 
$S^{*\left|n\right|} $. Let $p$ be the projection $\left|e_0\right\rangle \left\langle e_0\right|$. The number operator  
on $\ell^{2}\left(\bbn\right)$ is denoted by $N$ and defined  as $Ne_n :=ne_n$. The toeplitz algebra is denoted by $\scrt$.

\bdfn Let $A$ be a unital $C^*$-algebra. Then quantum double suspension of $A$ denoted by $\Sigma^{2} A$ is defined 
as the $C^*$-algebra generated by $A\otimes p$ and $1\otimes S$ in $A\otimes \scrt$.
\edfn
Let $\cla$ be a dense $*$-subalgebra of a $C^*$-algebra $A$. Define
\[
\Sigma^2(\cla)=\mbox{span}\left\{a\otimes k,1\otimes S^n : a\in \cla , k\in \cls(\ell^2(\bbn)), n \in \bbz \right\} 
\]
where $\cls(\ell^2(\bbn)) := \left\{(a_{mn}) : \sum_{m,n}(1+m+n)^p\left|a_{mn}\right| < \infty \quad \mbox{for} \quad p \in \bbn \right\}$.
Clearly $\Sigma^2(\cla)$ is a dense $*$ subalgebra of $\Sigma^2(A)$.
\bdfn Let $(\cla , \clh , D)$ be an odd spectral triple and denote the sign of $D$ by $F$. Then the spectral triple
$(\Sigma_{alg}^2(\cla), \clh \otimes \ell^2(\bbn), \Sigma^2(D):=((F\otimes 1)(\left|D\right|\otimes 1 + 1\otimes N))$ is 
called the quantum double suspension of the spectral triple $(\cla,\clh,D)$. For an even spectral triple  
$(\cla , \clh , D, \gamma)$ where $\gamma$ is the grading operator, the spectral triple
$(\Sigma_{alg}^2(\cla), \clh \otimes \ell^2(\bbn), \Sigma^2(D):=((F\otimes 1)(\left|D\right|\otimes 1 + 1\otimes N), 
\gamma \otimes 1)$ is called the quantum double suspension of the spectral triple $(\cla,\clh,D, \gamma)$.
\edfn
\brmrk  It should be emphasized that $\Sigma^2(D)$ can not be replaced by operators like $D\otimes I+I \otimes N$ for the simple reason that the latter operator do not have compact resolvent.
\ermrk
We use $D_0$ to denote the Dirac operator $\Sigma^2D$. We denote by  $d_0$ and $\delta_0$ the commutator with $D_0$  and 
$\left|D_0\right|$ respectively.  Given a multilinear functional $\phi$ defined  on $\cla$, we denote by 
  $\Sigma^2\phi$, the corresponding multilinear functional defined on $\Sigma^2\cla$. Now we proceed to our
main aim that is to describe local index formula for the quantum double suspension spectral triple in
terms of linear functionals appearing in local index formula for original spectral triple. 
To approach the problem,  
we need to put an extra condition namely WHKAE on the spectral triple.

Let $\phi : \left(0, \infty\right) \longrightarrow \bbc $ be a continuous function. We say that $\phi$ has an asymptotic expansion near $0$
if there exists a sequence of complex numbers $\left(a_r\right)_0^{\infty}$ such that given $N$ there exist $\epsilon, M$ such
that if $t \in \left(0, \infty \right)$
\begin{IEEEeqnarray}{rCl}
\left|\phi\left(t\right)-\sum_{r=0}^Na_rt^r\right|& \leq & Mt^{N+1} \nonumber
\end{IEEEeqnarray}
We write $\phi(t) \sim \sum_{r=0}^{\infty}a_rt^r$ as $t \rightarrow 0^{+}$. Note that the co-efficients $a_r$ are unique. For,
\begin{IEEEeqnarray}{rCl}
 a_N&=& \lim_{t\rightarrow 0^+} \frac{\phi\left(t\right)-\sum_{r=0}^N a_r t^r}{t^N}. \nonumber
\end{IEEEeqnarray}
If $\phi(t) \sim \sum_{r=0}^{\infty}a_rt^r$ as $t \rightarrow 0^{+}$  then $\phi$ can be extended continuously 
to $\left[0, \infty\right)$ simply by letting $\phi\left(0\right):=a_0$.

\bdfn \label{d1}
Let $\left( \cla. \clh, D\right)$ be a $p^+$-summable odd spectral triple for a $C^*$-algebra $A$ 
where $\cla$ is a dense $*$-subalgebra of $A$. We say that the spectral triple has the weak heat kernel 
asymptotic expansion property (WHKAE) of dimension $p$  if there exists a $*$-subalgebra $\clb \subset \cll ( \clh )$ such that
\begin{enumerate}
\item
$\clb$ contains $\cla$.
\item
The unbounded derivation $\delta :=\left[\left|D\right|,.\right]$ leaves $\clb$ invariant. Also the unbounded
derivation $d :=\left[D,.\right]$ maps $\cla$ into $\clb$.
\item
$\clb$ is invariant under the left multiplication by $F :=$sign$ D$.
\item
For every $b \in \clb$, the function $\tau_{p,b}(t) : \left[0, \infty\right) 
\longrightarrow \bbc$ defined by $\tau_{p,b}(t) = t^p Tr\left(be^{-t\left|D\right|}\right)$ has an asymptotic expansion.
\end{enumerate}
\edfn
\brmrk
\begin{enumerate}
\item
If the algebra $\cla$ is unital and the representation of $\cla$ on $\clh$ is unital then condition $3$ can be 
replaced by the condition $F \in \clb$.
\item
In case of an even spectral triple, we further demand $\clb$ to be invariant under left multiplication by the  
grading operator $\gamma$. If the algebra $\cla$ is unital and the representation of $\cla$ on $\clh$ is unital
then this condition can be replaced by the condition $\gamma \in \clb$.
\end{enumerate}
\ermrk 
It is known  that a spectral triple (odd or even) with 
WHKAE property is regular and has simple dimension spectrum.
\bppsn [Theorem 3.2, \cite{ChaSun-2011ab}]
Let $\left( \cla, \clh, D\right)$ be a $p^+$-summable odd spectral triple with 
WHKAE property of dimension $p$.
Then the spectral triple $\left( \cla. \clh, D\right)$ is regular and it has  finite simple dimension spectrum 
contained in $\left\{1,2,\cdots ,p\right\}$.
\eppsn 
\brmrk 
Above result holds for an even spectral triple. Similar proof will work. 
\ermrk
We state some results which relates the co-efficients of asymptotic expansion 
of the functions $t^p\trace\left(be^{-t\left|D\right|}\right)$ and 
$t^p\trace\left(be^{-t^2D^2}\right)$ with residues of the zeta functions $\trace(b\left|D\right|^{-2z})$ and  $\trace(b\left|D\right|^{-z})$.
\bppsn \label{res}
If $t^p\trace\left(be^{-t\left|D\right|}\right) \sim \sum_{r=0}^{\infty}b_rt^r$, then 
\begin{enumerate}
\item
$\mbox{Res}_{z=k}\trace\left(b\left|D\right|^{-z}\right)
=\frac{b_{p-k}}{\Gamma(k)}, \qquad  \mbox{for} \quad k \in \left\{1, 2,\cdots, p\right\}$.
\item
$\mbox{Res}_{z=k}\trace\left(b\left|D\right|^{-2z}\right)
=\frac{b_{p-2k}}{2\Gamma(2k)}, \qquad  \mbox{for} \quad k \in \left\{1/2, 2/2,\cdots, p/2\right\}$.
\item
$\mbox{Res}_{z=0}z^{-1}\trace\left(b\left|D\right|^{-z}\right)=\trace(b\left|D\right|^{-z})_{z=0}= b_p$.
\item
$\mbox{Res}_{z=0}z^{-1}\trace\left(b\left|D\right|^{-2z}\right)=\trace(b\left|D\right|^{-2z})_{z=0}= 2b_p$.
\end{enumerate}
\eppsn
\prf
It follows easily from remark $(3.3)$ \cite{ChaSun-2011ab}. \qed
\bppsn \label{compare}
If $t^p\trace\left(be^{-t\left|D\right|}\right) \sim \sum_{r=0}^{\infty}b_rt^r$ and 
$t^p\trace\left(be^{-t^2D^2}\right) \sim \sum_{r=0}^{\infty}b_r^{'}t^r$, then 
\[b_r=\frac{1}{\sqrt{\pi}}2^{p-r}\Gamma(\frac{p-r+1}{2})b_r^{'} \qquad 
\mbox{for}\quad r \in \{0,1,\cdots ,p\}.\]
\eppsn 

\prf For the proof, see proposition $3.5$ \cite{ChaSun-2011ab}.
\qed

The following proposition shows that WHKAE property is preserved under double suspension which ensures that 
the quantum double suspension spectral triple is regular and has simple dimension spectrum. 
We will consider odd spectral triple first. 

\bppsn [Theorem 4.5, \cite{ChaSun-2011ab}] \label{p3}
Let $(\cla , \clh , D)$ be an odd spectral triple with weak  heat kernel expansion property of dimension $p$.
Assume that the algebra $\cla$ is unital and the representation on $\clh$ is unital.
Then the spectral triple $(\Sigma_{alg}^2(\cla), \clh \otimes \ell^2(\bbn), \Sigma^2(D))$ also has 
the weak heat kernel expansion property with dimension $p+1$.
\eppsn
\prf We will give sketch of the proof. For more detail, see proposition $4.5$ \cite{ChaSun-2011ab}.
Let $\clb$ be a $*$ subalgebra of $\cll(\clh)$ for which $(1)-(4)$ of definition \ref{d1} holds. Define 
\[
\Sigma^2 \clb := \mbox{span}\left\{b \otimes k , 1\otimes S^n , F\otimes S^n : b\in \clb, k \in \cls(\ell^2(\bbn)), n \in \bbz \right\}
\]
where $\cls(\ell^2(\bbn)):=\left\{(a_{mn}): \sum_{m,n}(1+m+n)^p|a_{mn}| < \infty \right\}$.

Then $\Sigma^2\clb$ be a $*$ subalgebra of $\cll(\clh \otimes \ell^2(\bbn))$ for which $(1)-(4)$ of definition \ref{d1} 
satisfied for the spectral triple $(\Sigma_{alg}^2(\cla), \clh \otimes \ell^2(\bbn), \Sigma^2(D))$, taking  $\Sigma^2\clb$ as a bigger algebra.
\qed

 For $b \in \clb$ and $m \in \bbn$,  define
\[ 
\zeta_D^{(m)}(b)=\mbox{Res}_{z=m/2}\trace(b\left|D\right|^{-2z}).
\] 
Clearly $\zeta_D^{(m)}$ is a linear functional on $\clb$. The linear functionals $\psi_x^k$ and $\phi_n$ defined in the
previous section can be expressed in terms of the linear functionals $\left\{\zeta_D^{(m)}, m\in \bbn\right\}$. Hence
local index formula can also be written in these functionals. Therefore
to compute local index formula, it is enough to compute $\zeta_D^{(m)}$ on $\clb$ or certain subspace of $\clb$ defined as follows.
\bdfn
A subspace $\clb_{ess}$ of $\clb$ is called an essential subspace  if there exists another subspace $\clb_0$ with the following properties :
\begin{enumerate}
\item
$\clb=\clb_{ess} \oplus \clb_0$. 
%\item 
%$\zeta_D^{(m)}(b) =0$ for all $b \in \clb_0$ and $m \in \bbn$.
\item 
For all $b\in \clb_0$, $\trace(be^{-t\left|D\right|})$ has asymptotic power series expansion.
\end{enumerate}
\edfn
It follows from remark \ref{res} that $\zeta_D^{(m)}(b) =0$ for all $b \in \clb_0$ and $m \in \bbn$. Hence
it is enough to evaluate the linear functionals $\zeta_D^{(m)}$ on a essential subspace of $\clb$.
\bppsn \label{lres} 
Let $\clb_{ess}$ be an essential subspace of $\clb$. Then  
\[
\Sigma^2\clb_{ess} = \mbox{span}\left\{b\otimes k,1,F\otimes 1\quad \mbox{for} \quad  b\in \clb_{ess}, k \in \cls(\ell^2(\bbn)) \right\}.
\]
is an essential subspace of $\Sigma^2\clb$.
\eppsn

To prove  the proposition, we need to show that for $k \in \cls(\ell^2(\bbn))$,  $\trace(ke^{-tN})$ has
asymptotic expansion near $0$. Although exact expression for  asymptotic expansion of $\trace(ke^{-tN})$ is not
necessary here, we derive it for further use. For $k \in \cls(\ell^2(\bbn))$, define 
\[ 
\varphi_r(k) = \frac{(-1)^r}{r!}\sum_{i=0}^{\infty}k_{ii}i^r.
\]
where, $k_{ii}=\left\langle ke_i,e_i\right\rangle$ and $\left\{e_i\right\}_{i=0}^{\infty}$ is the standard
orthonormal basis of $\ell^2(\bbn)$. Observe that the sum given above converges as $k \in \cls(\ell^2(\bbn))$
and hence, $\varphi_r(k)$ is  well- defined linear functional on $\cls(\ell^2(\bbn))$.
\bppsn \label{asy}
For $k \in \cls(\ell^2(\bbn))$,
\[
\trace(ke^{-tN}) \sim \sum_{r=0}^{\infty}\varphi_r(k)t^r.
\]
\eppsn 
\prf
It is easy to verify the expression for  $k = \left|e_i\rangle \langle e_j\right|$. By linearity, it holds for all finite rank operators. 
For $k \in \cls(\ell^2(\bbn))$, the finite rank operators $k_n:= \sum_{i=0}^n k_{ij}\left|e_i\rangle \langle e_j\right|$ converges to 
$k$ as $n\rightarrow \infty$. Thus $\lim_{n\rightarrow \infty}\phi_r(k_n)= \phi_r(k)$.
\qed

\prf (of proposition \ref{lres}): Let $\clb_0$ be 
a complementary subspace of $\clb_{ess}$ such that $\trace(be^{-t\left|D\right|}) =0$ for all $b \in \clb_0$. Define 
\[
\Sigma^2\clb_0=\mbox{span}\left\{b\otimes k,1\otimes S^n ,F\otimes S^n\quad \mbox{for} \quad  b\in \clb_0, k \in \cls(\ell^2(\bbn)),
n\in \bbz-\left\{0\right\} \right\}
\]
Clearly $\Sigma^2\clb=\Sigma^2\clb_{ess} \oplus \Sigma^2\clb_0$. 
 For $n \in \bbz-\left\{0\right\}$, $\trace((1\otimes S^n)e^{-t\left|D_0\right|})=0$ and $\trace((F\otimes S^n)e^{-t\left|D_0\right|})=0$.
%\begin{IEEEeqnarray}{rCl}
%\trace(1\otimes S^n\left|D_0\right|^{-2z})&=&0. \nonumber \\
%\trace(F\otimes S^n\left|D_0\right|^{-2z})&=&0. \nonumber
%\end{IEEEeqnarray}
Now for $b \in \clb_0$, it follows from remark \ref{res} that  $\trace(be^{-t\left|D\right|})$ has an asymptotic expansion near $0$.
Also, from proposition \ref{asy}, it follows  that $\trace(ke^{-t\left|N\right|})$ has  an asymptotic expansion near $0$.
Therefore, $\trace((b\otimes k)e^{-t\left|D_0\right|}) =\trace(be^{-t\left|D\right|})\trace(ke^{-t\left|N\right|})$ has 
asymptotic expansion near $0$. This completes the proof. \qed\\
Expansion of the  linear functionals  $\zeta_{D_0}^{(m)}$ and $\Sigma^2\phi_x^{(m)}$   
involves many residues of the form $\zeta_{D}^{(s)}$. So, removing unnecessary terms using structure of the algebra $\clb$
will be very crucial for avoiding difficult calculations. Along this direction, let
$I_0 \subset I_1 \subset \cdots  \subset I_p=\clb$ be ideals in $\clb$ such that 
\begin{enumerate} 
\item
$a \in I_{\ell} \Longrightarrow da \in I_{\ell}\quad \mbox{and} \quad  \delta(a) \in I_{\ell}$ ,
\item
for every $a \in I_{\ell}$, $t^{\ell}\trace(ae^{-t\left|D\right|})$ has asymptotic expansion near $0$. 
\end{enumerate} 
Set $J_{\ell}=I_{\ell}\otimes \cls(\ell^2(\bbn)), \quad \mbox{for} \quad \ell \leq p \quad$ and
$\quad J_{p+1}=\Sigma^2 \clb$. Then  $J_0 \subset J_1 \subset \cdots , J_p \subset J_{p+1}$ are
ideals in $\Sigma^2\clb$ with the following property.  
\begin{enumerate}
\item
$a \in J_{\ell} \Longrightarrow d_0a \in J_{\ell} \quad \mbox{and} \quad \delta_0(a) \in J_{\ell}$,
\item
for every $a \in J_{\ell}$,  $t^{\ell}\trace(ae^{-t\left|D_0\right|})$ has asymptotic expansion near $0$. 
\end{enumerate} 
Note that for $b \in I_{\ell}$, $\mbox{Res}_{z=\frac{k}{2}}\trace(b\left|D\right|^{-2z})=0$ and 
for $\tilde{b} \in J_{\ell}$, $\mbox{Res}_{z=\frac{k}{2}}\trace(\tilde{b}\left|D_0\right|^{-2z})=0$ for $k>\ell$.

\bppsn \label{l1}
For $b \in I_m, k \in \cls(\ell^2(\bbn))$ and $s \in \left\{1,2,\cdots, m\right\}$, one has
\[ \zeta_{D_0}^{(s)}(b\otimes k)=\frac{1}{\Gamma(s)}\sum_{y=0}^{m-s}\Gamma(s+y)\zeta_D^{(s+y)}(b)\varphi_y(k).
\]
\eppsn
\prf
Let $t^m\trace(be^{-t\left|D\right|}) \sim \sum_{r=0}^{\infty}b_rt^r$. Then
\begin{IEEEeqnarray}{rCl}
t^m\trace((b\otimes k)e^{-t\left|D_0\right|}) &=&t^m\trace(be^{-t\left|D\right|})\trace(ke^{-tN}).\nonumber \\
&\sim & \sum_{r=0}^{\infty}b_rt^r \sum_{r=0}^{\infty}\varphi_r(k)t^r. \nonumber \\
&\sim & \sum_{r=0}^{\infty}c_rt^r \nonumber.
\end{IEEEeqnarray} 
where $ c_r= \sum_{y=0}^r\varphi_y(k)b_{r-y}$. Hence for $s\in \left\{1,2,\cdots,m\right\}$, we get 
\begin{IEEEeqnarray}{rCl}
\zeta_{D_0}^{(s)}(b\otimes k)&=&\mbox{Res}_{z=\frac{s}{2}}\trace((b\otimes k)\left|D_0\right|^{-2z}). \nonumber \\
&=&\frac{1}{2\Gamma(r)}c_{m-s}. \nonumber \\
&=&\frac{1}{2\Gamma(s)}\sum_{y=0}^{m-s} \varphi_y(k)b_{m-(s+y)}. \nonumber \\
&=&\frac{1}{\Gamma(s)}\sum_{y=0}^{m-s}\Gamma(s+y)\zeta_D^{(s+y)}(b)\varphi_y(k). \nonumber 
\end{IEEEeqnarray} 
\qed  \\
Observe that number of terms in the expression of $\zeta_{D_0}^{(s)}(b\otimes k)$ is depending on $m$. Let 
\begin{align} \label{chap8-eqn-asy}
 t^p\trace(e^{-t\left|D\right|}) &\sim  \sum_{r=0}^{\infty}u_rt^r, 
& t^p\trace(Fe^{-t\left|D\right|}) &\sim  \sum_{r=0}^{\infty}v_rt^r, 
& t\trace(e^{-tN})  &\sim  \sum_{r=0}^{\infty}n_rt^r. 
\end{align}

\bppsn \label{l2}
For $s \in \left\{2,3,\cdots,p+1\right\}$, one has
\[
 \zeta_{D_0}^{(s)}(1) = \frac{1}{\Gamma(s)}\sum_{y=0}^{p+1-s}\Gamma(s+y)\zeta_D^{(s+y)}(1)n_y, \quad \quad
\zeta_{D_0}^{(s)}(F\otimes 1)=\frac{1}{\Gamma(s)}\sum_{y=0}^{p+1-s}\Gamma(s+y)\zeta_D^{(s+y)}(F)n_y. 
\]
\eppsn
\prf
We have
\begin{IEEEeqnarray}{rCl}
t^{p+1}\trace(e^{-t\left|D_0\right|}) &=&t^p\trace(e^{-t\left|D\right|})\trace(e^{-tN}).\nonumber \\
&\sim & \sum_{r=0}^{\infty}u_rt^r \sum_{r=0}^{\infty}n_rt^r. \nonumber \\
&\sim & \sum_{r=0}^{\infty}\tilde{u}_rt^r \nonumber.
\end{IEEEeqnarray}
where $ \tilde{u}_r= \sum_{y=0}^rn_{y}b_{r-y}$. Hence for $s \in \left\{2,3,\cdots,p+1\right\}$, we have
\begin{IEEEeqnarray}{rCl}
\zeta_{D_0}^{(s)}(1)&=&\mbox{Res}_{z=\frac{s}{2}}\trace(\left|D_0\right|^{-2z}).\nonumber \\
&=&\frac{1}{2\Gamma(s)}\tilde{u}_{p+1-s}. \nonumber \\
&=&\frac{1}{2\Gamma(s)}\sum_{y=0}^{p+1-s} n_yu_{p+1-(s+y)}. \nonumber \\
&=&\frac{1}{\Gamma(s)}\sum_{y=0}^{p+1-s}\Gamma(s+y)\zeta_D^{(s+y)}(1)n_y. \nonumber 
\end{IEEEeqnarray}
Similar calculation will prove that
\begin{IEEEeqnarray}{rCl}
\zeta_{D_0}^{(s)}(F\otimes 1) &=&\frac{1}{\Gamma(s)}\sum_{y=0}^{p+1-s}\Gamma(s+y)\zeta_D^{(s+y)}(F)n_y. \nonumber 
\end{IEEEeqnarray} \qed 
\bppsn \label{l3} Let $u_p$, $v_p$ and $n_i$'s be given by equation $(\ref{chap8-eqn-asy})$. Then one has
\[
 \zeta_{D_0}^{(1)}(1)=n_0u_p+ \sum_{y=1}^{p}\Gamma(y)\zeta_D^{(y)}(1)n_y, \quad \quad
\zeta_{D_0}^{(1)}(F\otimes 1)=n_0v_p+\sum_{y=1}^{p}\Gamma(y)\zeta_D^{(y)}(F)n_y. 
\]
\eppsn 
\prf 
Let  $\tilde{u}_{p+1}$ is as in the  previous proposition. Then
\begin{IEEEeqnarray}{rCl}
\zeta_{D_0}^{(1)}(1)&=&\mbox{Res}_{z=\frac{1}{2}}\trace(\left|D_0\right|^{-2z}).\nonumber \\
&=&\frac{1}{2\Gamma(1)}\tilde{u}_{p+1}.  \nonumber \\
&=&\frac{1}{2}\sum_{y=0}^{p} n_yu_{p+1-y}. \nonumber \\
&=& n_0u_p+ \sum_{y=1}^{p}\Gamma(y)\zeta_D^{(y)}(1)n_y. \nonumber 
\end{IEEEeqnarray}
Similar calculations will prove the other part of the assertion. \qed

 The above calculations shows that the linear functionals $\left\{\zeta_{D_0}^{(k)}\right\}_{k=1}^{p+1}$ on 
 $\Sigma^2\clb$ can be computed by the linear functionals $\left\{\zeta_{D}^{(k)}\right\}_{k=1}^{p}$ on $\clb$ and 
 two values $\zeta_{D_0}^{(1)}(1)$ and $\zeta_{D_0}^{(1)}(F\otimes 1)$ (or $u_p$ and $v_p$).  We will now put
 these data altogether and derive exact expression of local index formula for the spectral 
 triple $(\Sigma^2\clb,\clh \otimes \ell^2(\bbn),D_0) $. For $x \in \bbn^n$ and $\ell \in \bbn$, define
\begin{IEEEeqnarray}{rCl}
\Sigma^2\psi_{x}^{\ell}(\tilde{b_0},\tilde{b_1}, \cdots ,\tilde{b_n})
&=&\mbox{Res}_{z=\left(n+\left|x\right|+\ell \right)/2}\trace \left(\tilde{b_0}\delta_{D_0}^{x_1}(\tilde{b_1})
\cdots \delta_{D_0}^{x_n}(\tilde{b_n})\left|D_0\right|^{-2z}\right) \nonumber 
\end{IEEEeqnarray}

where $\tilde{b}_i \in \Sigma^2\clb$. To evaluate this, it is enough to take elements of the
type $a\otimes k$ and $1\otimes S^n$ which we commonly write in the form $a\otimes c$. Let
$a^{(1)}=da,a^{(2)}=Fa, a^{(3)}=aF,c^{(1)}=c,c^{(2)}=Nc$ and $c^{(3)}=-cN$. Then 
\begin{IEEEeqnarray}{lCl}
d_0(a\otimes c)&=& da\otimes c + aF\otimes cN - Fa\otimes Nc =a^{(1)}\otimes c^{(1)} + 
a^{(2)}\otimes c^{(2)} + a^{(3)}\otimes c^{(3)}\nonumber 
\end{IEEEeqnarray}

\blmma \label{lif1}
Let $x \in \bbn^n$. If for some $i \in \left\{0,1,\cdots ,n\right\}$ and for some $m \in 
\left\{1,2,\cdots,p\right\}$, $a_i\otimes c_i \in J_m$
then 
\begin{IEEEeqnarray}{rCl}
\Sigma^2\psi_{x}^{0}(a_0\otimes c_0,d_0(a_1\otimes c_1), \cdots ,d_0(a_n\otimes c_n))
&=&\sum_{\substack{1\leq j_1\leq 3 \\ \cdots \\ 1\leq j_n\leq 3}} \sum_{\substack{0 \leq r_1 \leq x_1 \\ \cdots \\  0 \leq r_n \leq x_n}}\sum_{k=0}^{m-(n+\left|x\right|)} \prod_{i=1}^n{x_i\choose r_i}
\frac{\Gamma(n+\left|x\right|+k)}{\Gamma(n+\left|x\right|)} \nonumber \\
& &\psi_{r}^{(k+\left|x\right|-\left|r\right|)}(a_0,a_1^{(j_1)},\cdots,a_n^{(j_n)})\varphi_j(c_0\prod_{i=1}^n
\delta_N^{x_i-r_i}(c_i^{(j_i)})).\nonumber
\end{IEEEeqnarray}
\elmma
\prf
%\textbf{Case 1}: $a_i\otimes c_i \in J_m, m\leq p$ for some $i \in \left\{0,1,\cdots,n\right\}$. \\
 It is easy to see that $\delta_{0}^{n}(a\otimes c)= \sum_{r=0}^n {n\choose r}\delta^r(a)\otimes \delta_N^{n-r}(c)$.
 Using this, we get
\begin{IEEEeqnarray}{lCl}
\Sigma^2\psi_{x}^{0}(a_0\otimes c_0,d_0(a_1\otimes c_1), \cdots ,d_0(a_n\otimes c_n))\nonumber \\ 
=\mbox{Res}_{z=\left(n+\left|x\right| \right)/2}\trace\left(a_0\otimes c_0\prod_{i=0}^n\delta_{0}^{x_i}d_0
(a_i\otimes c_i)
\left|D_0\right|^{-2z}\right) \nonumber \\
=\sum_{\substack{1\leq j_1\leq 3 \\ \cdots \\ 1\leq j_n\leq 3}}\mbox{Res}_{z=\left(n+\left|x\right| \right)/2}\trace\left(a_0\otimes c_0\prod_{i=0}^n\delta_{0}^{x_i} (a_i\otimes c_i)\left|D_0\right|^{-2z}\right) \nonumber \\
=\sum_{\substack{1\leq j_1\leq 3 \\ \cdots \\ 1\leq j_n\leq 3}} \sum_{\substack{0 \leq r_1 \leq x_1 \\ \cdots \\  0 \leq r_n \leq x_n}}\prod_{i=1}^n {x_i\choose r_i}\mbox{Res}_{z=\left(n+\left|x\right| \right)/2}
\trace\left(a_0\prod_{i=1}^n \delta^{r_i}(a_i^{(j_i)})\otimes c_0\prod_{i=1}^n \delta_N^{x_i-r_i}(c_i^{(j_i)})
\left|D_0\right|^{-2z}\right) \nonumber \\
=\sum_{\substack{1\leq j_1\leq 3 \\ \cdots \\ 1\leq j_n\leq 3}} \sum_{\substack{0 \leq r_1 \leq x_1 \\ \cdots \\  0 \leq r_n \leq x_n}}\prod_{i=1}^n {x_i\choose r_i} 
\zeta_{D_0}^{(n+\left|x\right|)}\left(a_0\prod_{i=1}^n \delta^{r_i}(a_i^{(j_i)})\otimes c_0\prod_{i=1}^n
\delta_N^{x_i-r_i}(c_i^{(j_i)})\right)\nonumber \\ 
=\sum_{\substack{1\leq j_1\leq 3 \\ \cdots \\ 1\leq j_n\leq 3}} \sum_{\substack{0 \leq r_1 \leq x_1 \\ \cdots \\  0 \leq r_n \leq x_n}}
\sum_{k=0}^{m-(n+\left|x\right|)}\frac{\Gamma(n+\left|x\right|+k)}{\Gamma(n+\left|x\right|)} \prod_{i=1}^n{x_i\choose r_i}
\psi_{r}^{(k+\left|x\right|-\left|r\right|)}(a_0,a_1^{(j_1)},\cdots,a_n^{(j_n)})\varphi_k(c_0\prod_{i=1}^n 
\delta_N^{x_i-r_i}(c_i^{(j_i)})).\nonumber
\end{IEEEeqnarray} \qed

\bppsn \label{lif2}
Let $x \in \bbn^n$. For $m_0, m_1, \cdots ,m_n \in \bbz$ such that $ \sum_{i=0}^nm_i= 0$, one has
\begin{IEEEeqnarray}{rCl}
\Sigma^2\psi_{x}^{0}(1\otimes S^{m_0},d_0(1\otimes S^{m_1}), \cdots ,d_0(1\otimes S^{m_n}))&=&-\prod_{i=1}^nm_i^{x_i+1}\zeta_{D_0}^{(n+\left|x\right|)}(F\otimes 1). \nonumber 
\end{IEEEeqnarray}
\eppsn
\prf
 Observe that $d_0(1\otimes S^m)=-mF\otimes S^m$ and $\quad \delta_0^{x}d_0(1\otimes S^m)=-m^{x+1}F\otimes S^m$.
Hence we have
\begin{IEEEeqnarray}{lCl}
\Sigma^2\psi_{x}^{0}(1\otimes S^{m_0},d_0(1\otimes S^{m_1}), \cdots ,d_0(1\otimes S^{m_n})) \nonumber \\
=\mbox{Res}_{z=\left(n+\left|x\right| \right)/2}\trace\left(1\otimes S^{m_0}
\prod_{i=0}^n\delta_{0}^{x_i}d_0(1\otimes S^{n_i})\left|D_0\right|^{-2z}\right) \nonumber \\
=(-1)^n\prod_{i=1}^nm_i^{x_i+1}\mbox{Res}_{z=\left(n+\left|x\right| \right)/2}
\trace\left((F^n\otimes S^{\sum_{i=0}^nm_i})\left|D_0\right|^{-2z}\right)\nonumber \\
=-\prod_{i=1}^nm_i^{x_i+1}\mbox{Res}_{z=\left(n+\left|x\right| \right)/2}
\trace\left((F\otimes 1)\left|D_0\right|^{-2z}\right).\nonumber \\
=-\prod_{i=1}^nm_i^{x_i+1}\zeta_{D_0}^{(n+\left|x\right|)}(F\otimes 1)  \nonumber
\end{IEEEeqnarray} 
Note that since we are taking odd spectral triple, $n$ is odd in this case. 
\qed

It is easy to see that $ \sum_{i=0}^nm_i \neq 0$,
then $\Sigma^2\psi_{x}^{0}(1\otimes S^{m_0},d_0(1\otimes S^{m_1}), \cdots ,d_0(1\otimes S^{m_n}))=0$. 
Putting these results in equation $(\ref{e4})$, we get the local index formula for the 
spectral triple $(\Sigma^2\cla, \clh \otimes \ell^2(\bbn),\Sigma^2D)$.

Now we will consider even spectral triple $(\cla,\clh,D,\gamma)$ and its quantum double
suspension \\$(\Sigma^2\cla,\clh \otimes \ell^2(\bbn),\Sigma^2D,\gamma \otimes 1)$. In this case, we define
\[
\Sigma^2\clb:=\mbox{span} \left\{b\otimes k, 1,F\otimes S^n, \gamma \otimes S^n, \gamma F\otimes
S^n : b \in \clb, k \in \cls(\ell^2(\bbn)),n \in \bbz \right\}.
\]
where $\clb$ is a *-subalgebra of $\cll(\clh)$ for which all conditions for WHKAE property of
 spectral triple $(\cla,\clh,D,\gamma)$ holds. Since, $F \otimes S^n$ and $\gamma F\otimes S^n$ are  odd operators, 
 $\zeta_{D_0}^{(m)}(F \otimes S^n)=0$ and $\zeta_{D_0}^{(m)}(\gamma F \otimes S^n)=0$ for all $m \in 
 \bbn$ and $n \in \bbz$. 
 Using this, one can easily check that if $\clb_{ess}$ is an essential subspace of
 $\clb$ then $\Sigma^2\clb_{ess}$ defined below 
 will be an essential subspace of $\Sigma^2\clb$.
\[
\Sigma^2\clb_{ess}= \mbox{span} \left\{b\otimes k, 1,
\gamma \otimes 1 : b \in \clb_{ess}, k \in \cls(\ell^2(\bbn))\right\}.
\]

Let 
\begin{IEEEeqnarray}{rCl} \label{chap8-eqn-gamma} 
 t^p\trace(\gamma e^{-t\left|D\right|}) \sim  \sum_{r=0}^{\infty}w_rt^r. 
\end{IEEEeqnarray}
Then lemma \ref{l1}, \ref{l2} and  \ref{l3} will follow with $F$ and $v_r$ replaced by 
$\gamma$ and $w_r$ respectively. For local index formula, we state the following results. Since proof of these results 
are similar to that in odd case, we omit it.
\blmma \label{lif3} 
Let $x \in \bbn^n$ where $n$ is positive even integer. If for some $i \in \left\{0,1,\cdots ,n\right\}$ 
and for some $m \in \left\{1,2,\cdots,p\right\}$, $a_i\otimes c_i \in J_m$  then 
\begin{IEEEeqnarray}{lCl}
\Sigma^2\psi_{x}^{0}((\gamma \otimes 1)(a_0\otimes c_0),d_0(a_1\otimes c_1), \cdots ,d_0(a_n\otimes c_n)) \nonumber \\
=\sum_{\substack{1\leq j_1\leq 3 \\ \cdots \\ 1\leq j_n\leq 3}} \sum_{\substack{0 \leq r_1 \leq x_1 \\ \cdots \\  0 \leq r_n \leq x_n}}\sum_{k=0}^{m-(n+\left|x\right|)}\prod_{i=1}^n{x_i\choose r_i}
\frac{\Gamma(n+\left|x\right|+k)}{\Gamma(n+\left|x\right|)} \nonumber \\
\psi_{r}^{(k+\left|x\right|-\left|r\right|)}(\gamma a_0,a_1^{(j_1)},\cdots,a_n^{(j_n)})
\varphi_k(c_0\prod_{i=1}^n \delta_N^{x_i-r_i}(c_i^{(j_i)})).\nonumber
\end{IEEEeqnarray}
\elmma
\bppsn \label{lif4}
Let $x \in \bbn^n$ where $n$ is positive even integer. For $m_0, m_1, \cdots ,m_n \in \bbz$ such that $ \sum_{i=0}^nm_i= 0$, 
\begin{IEEEeqnarray}{rCl}
\Sigma^2\psi_{x}^{0}((\gamma \otimes 1)(1\otimes S^{m_0}),d_0(1\otimes S^{m_1}),
\cdots ,d_0(1\otimes S^{m_n}))&=&-\prod_{i=1}^nm_i^{x_i+1}\zeta_{D_0}^{(n+\left|x\right|)}(\gamma \otimes 1). \nonumber 
\end{IEEEeqnarray}
\eppsn 
For $\sum_{i=0}^nm_i \neq 0$, it is easy to check that 
$\Sigma^2\psi_{x}^{0}((\gamma \otimes 1) (1\otimes S^{m_0}),d_0(1\otimes S^{m_1}), \cdots ,d_0(1\otimes S^{m_n}))=0$.
Putting these results in equation $(\ref{e5})$, we can compute the functionals $\Sigma^2\phi_n$ for all $n >0$ and even.
Moreover, we need to calculate $\Sigma^2\phi_0(\tilde{b}):=\mbox{Res}_{z=0}z^{-1}\trace(\tilde{b}
\left|D_0\right|^{-2z})$ where $\tilde{b} \in \Sigma^2\clb$.
\bppsn \label{l6}
For $b \in I_m$ and $s \in \left\{1,2,\cdots, m\right\}$, one has
\[ \Sigma^2\phi_0(b\otimes k)=\phi_0(b)\varphi_0(k)+\frac{1}{\Gamma(s)}\sum_{y=1}^{m}
\Gamma(y)\zeta_D^{(y)}(b)\varphi_y(k).
\]
\eppsn
\bppsn \label{l5} Let $w_p$ and $w_{p+1}$ be given by equation $(\ref{chap8-eqn-gamma})$. Assume that  
$u_p$, $u_{p+1}$ and $n_i$'s are given by $(\ref{chap8-eqn-asy})$. Then one has
\begin{IEEEeqnarray}{rCl}
\Sigma^2\phi_0(1)&=&2n_0u_{p+1}+2n_1u_p+ \sum_{y=1}^{p}\Gamma(y)\zeta_D^{(y)}(1)n_y. \nonumber \\
\Sigma^2\phi_0(\gamma \otimes 1)&=&2n_0v_{p+1}+2n_1v_p+\sum_{y=1}^{p}\Gamma(y)\zeta_D^{(y)}(\gamma)n_y. \nonumber 
\end{IEEEeqnarray} 
\eppsn
Putting these results in equation $(\ref{e5})$, we get local index formula for the even spectral 
triple  $(\Sigma_{alg}^2(\cla), \clh \otimes \ell^2(\bbn), \Sigma^2(D), \gamma \otimes 1)$.

\brmrk 
 We computed linear functionals for finite linear combination of elements from the set
 $\left\{b\otimes \cls(\ell^2(\bbn)), 1\otimes S^n, F\otimes S^n: b\in B, n\in \bbz \right\}$. We can 
 extend these linear functionals to infinite linear combination of these element in which co-efficients are
 rapidly decreasing. Hence we can apply these results to the spectral triple with topological weak heat kernel
 expansion property defined as in \cite{ChaSun-2011ab} . 
\ermrk

\section{SOME EXAMPLES}

\subsection{Local index formula for $\Sigma^2C(\cls^2)$}
 Consider the $2$-dimensional sphere $\cls^2$ with usual orientation. The usual sphereical coordinates on $\cls^2$ are:
\[ 
p=(\sin\theta \cos\phi, \sin\theta \sin\phi, \cos\theta) \in \cls^2.
\] 
The poles are $N=(0,0,1)$ and $S=(0,0,-1)$. Let $U_N = \cls^2-\left\{N\right\}$ and $U_S = \cls^2-\left\{S\right\}$ be the 
two charts of $\cls^2$. Consider the  stereographic projections  $p \mapsto z$ from $U_N$ to  $\bbc$ and   
$p \mapsto \xi$ from $U_S$ to  $\bbc$ given by
\[ 
z:=e^{-i\phi}\cot\frac{\theta}{2}, \qquad  \xi:=e^{i\phi}\tan\frac{\theta}{2},
\]
such that $z=\frac{1}{\xi}$ on $U_N\cap U_S$. Let $L^+$ be the tautological 
line bundle coming from $\cls^2\cong \bbc P^1$. Let $L^-$ be the dual line bundle of $L^+$. 
One can show that  $L^+$ and $L^-$ are two nonisomorphic nontrivial complex line bundles. 
Define $L=L^+ \oplus L^-$. Denote by $\Gamma(\cls^2, L)$,  $\Gamma(\cls^2, L^+)$ and $\Gamma(\cls^2, L^-)$ the sets
of smooth sections of the bundles $L$ ,  $L^+$ and $L^-$ respectively.  A spinor  on $\cls^2$ (smooth section of
the bundle $L$) is given by  two pairs of smooth functions
\[
\left({\begin{array}{cc}
\psi_N^+(z,\bar{z})\\
\psi_N^-(z,\bar{z})\\
\end{array}}\right)
\mbox{on} \quad U_N, \qquad 
\left({\begin{array}{cc}
\psi_S^+(\xi,\bar{\xi})\\
\psi_S^-(\xi,\bar{\xi})\\
\end{array}}\right)
\mbox{on} \quad U_S
\]
satisfying the following properties:
 
\begin{enumerate}
\item
$\psi_N^+(z,\bar{z}) =(\frac{\bar{z}}{z})^{\frac{1}{2}}\psi_S^+(z^{-1},\bar{z}^{-1}) \quad$ for $ z \neq 0$.
\item 
$\psi_N^-(z,\bar{z}) =(\frac{z}{\bar{z}})^{\frac{1}{2}}\psi_S^-(z^{-1},\bar{z}^{-1}) \quad$ for $ z \neq 0$. 
\item
$\psi_N^+$ and $\psi_N^-$  are regular at $ z=0$.
\item
$\psi_S^+$ and $\psi_S^-$  are regular at $ \xi=0$.
\end{enumerate}
Note that $\psi_N^+$ on $U_N$ and $\psi_S^+$ on $U_S$ gives a smooth section of the bundle $L^+$ and $\psi_N^-$ on $U_N$ and 
$\psi_S^-$ on $U_S$ gives a smooth section of the bundle $L^-$. One can show that $\Gamma(\cls^2, L)=\Gamma(\cls^2, L^+)\oplus \Gamma(\cls^2, L^-)$.
The scalar product on $\Gamma(\cls^2, L^+)$ is given by 
\[
\left\langle \phi, \psi \right\rangle = \int_{\cls^2}\left\langle \phi(p), \psi(p)\right\rangle_{p}\nu_g,
\]
where $\left\langle . ,. \right\rangle_p$ is the standard scalar product on $\bbc^2$ and $\nu_g$ is the
Riemannian volume form on $\cls^2$. On completion in the norm $\left\|\phi\right\|=\sqrt{\left\langle \phi,\phi\right\rangle}$, 
we get the Hilbert space $\clh^+:=L^2(\cls^2,L^+)$, the $L^2$-spinor of the bundle $L^+$. In a similar way one can construct
the Hilbert space $\clh^-:=L^2(\cls^2,L^-)$, the $L^2$-spinor of the bundle $L^-$ and  $\clh:=L^2(\cls^2,L)$,  the $L^2$-spinor 
of the bundle $L$. It is easy to see that $\clh=\clh^+\oplus \clh^-$. Let $\gamma$ be the grading on $\clh$ associated with this 
decomposition i.e.\  $\gamma(h)=h$ if $h \in \clh^+$ and  
$\gamma(h)=-h$ if $h \in \clh^-$. We represent $C^{\infty}(\cls^2)$ on $\clh$ 
by multiplication operators.  Let $D$ be the Dirac operator associated 
with  the  Levi-Civita connection. It is given over $U_N$ by
\[
-i
\left({\begin{array}{cc}
0, & (1+z\bar{z})\frac{\partial}{\partial z}-\frac{1}{2}\bar{z} \\
(1+z\bar{z})\frac{\partial}{\partial \bar{z}}-\frac{1}{2}z, & 0\\
\end{array}}\right).
\] 
A similar expression is valid over $U_S$ by replacing $z$ and $\bar{z}$ by $\xi$ and $\bar{\xi}$ respectively
and by changing overall $(-i)$ factor to $i$.
Then  $(C^{\infty}(\cls^2), \clh , D, \gamma)$ is the classical  even spectral triple of $C(\cls^2)$. For a complete description,   
we refer the reader to Varilly (\cite{Var-2006aa}, page 98-102).
Chakraborty and Sundar (\cite{ChaSun-2011ab},  page 15) showed that this spectral triple has the WHKAE property with dimension $2$ and hence
it is regular with dimension spectrum  contained in $\left\{1,2\right\}$. It follows from Proposition \ref{p3} that 
the quantum double suspension spectral triple $(\Sigma^2(C^{\infty}(\cls^2)), \clh \otimes \ell^2(\bbn), D_0:=\Sigma^2D),\gamma \otimes 1)$
also has the WHKAE property and its dimension spectrum is contained in $\left\{1,2,3\right\}$. 
Our aim is to give an explicit description of the
local index formula for this spectral triple. We will use symbol calculus as our main tool (see \cite{Gil-1995aa}, \cite{Var-2006aa}).

Since the  dimension spectrum is contained in $\left\{1,2,3\right\}$, it 
follows from equation 
$(\ref{e5})$ that $\Sigma^2\phi_{2n}=0$ for $n>1$.
So, we  need to compute $\Sigma^2\phi_0$ and $\Sigma^2\phi_2$ only. We first find the asymptotic expansion 
of $t^2\trace(\gamma a e^{-t^2D^2})$ for all
$a \in C^{\infty}(\cls^2)$.
\bppsn \label{sym}
 For all $a \in C^{\infty}(\cls^2)$, one has
\begin{IEEEeqnarray}{lCl}
t^2\trace(\gamma a e^{-t^2D^2}) & \sim & 0. \nonumber
%t^2\trace(\gamma a e^{-t\left|D\right|}) & \sim & 0. \nonumber 
\end{IEEEeqnarray}
\eppsn
\prf
 The symbol attached 
with the Dirac operator $D$  over a local chart is 
\begin{displaymath}
 p^D(x, \xi)=\begin{pmatrix}
   0 & \xi \\
   \bar{\xi} & 0 
  \end{pmatrix}  + 
  \frac{i}{2} \begin{pmatrix}
   0 & x_1-ix_2 \\
   x_1+ix_2 & 0 
  \end{pmatrix}.
\end{displaymath}
 Hence symbol for the pseudodifferential operator $D^2$ over the same chart will be given by
	\[
	p^{D^2}(x,\xi)=\begin{pmatrix}
   \left|\xi\right|^2  & 0\\
   0 & \left|\xi\right|^2  \\
  \end{pmatrix} + 
  \frac{i}{2}\begin{pmatrix}
   \xi_1x_1+\xi_2x_2 & 0 \\
   0 & \xi_1x_1+\xi_2x_2
  \end{pmatrix}  + \begin{pmatrix}
   1-\left|x\right|^2/4 & 0 \\
   0 & 1-\left|x\right|^2/4
  \end{pmatrix}. 
	\]
Let $p_k(x,\xi)$ be homogeneous polynomials in variable $\xi$ of degree $k$ such that $p^{D^2}(x,\xi)=\sum_{k=0}^2p_k(x,\xi)$. Then
	\[
	p_0(x,\xi)=\begin{pmatrix}
   1-\left|x\right|^2/4 & 0 \\
   0 & 1-\left|x\right|^2/4
  \end{pmatrix} ,\quad
   p_1(x,\xi)=\frac{i}{2}\begin{pmatrix}
   \xi_1x_1+\xi_2x_2 & 0 \\
   0 & \xi_1x_1+\xi_2x_2
  \end{pmatrix},
	\]
	and the leading symbol is 
	\[
	 p_2(x,\xi)=\begin{pmatrix}
   \left|\xi\right|^2  & 0\\
   0 & \left|\xi\right|^2  
  \end{pmatrix}.
\]
 Note that the leading symbol $p_2(x, \psi)$ is positive definite and is a scalar matrix.
 Let $K(t,x,y)$ be the kernel of $e^{-tD^2}$.
 Using Lemma $1.7.4$ in \cite{Gil-1995aa}, we get
\[
K(t,x,x) \sim \sum_{n=0}^{\infty}t^\frac{n-2}{2}e_n(x) \quad \mbox{as}\quad t\rightarrow 0^+.
\] 
Exact expression for $e_n(x)$ is given in (\cite{Gil-1995aa} page 54). 
Since  $p_k(x, \xi)$ is a scalar matrix for all $k \in \left\{0,1,2\right\}$,  
$e_n(x)$ will be a scalar matrix. It follows from (Lemma $1.7.7$, \cite{Gil-1995aa})
 that the kernel of  $\gamma ae^{-tD^2}$ is $\gamma aK(t,x,y)$. Hence
\begin{IEEEeqnarray}{rCl}
t\trace(\gamma ae^{-tD^2})&=&t\int_{\cls^2}\trace( \gamma a(x)K(t,x,x))\nu_g \nonumber \\
&\sim&t\sum_{n=0}^{\infty}t^{\frac{n-2}{2}}\int_{\cls^2}\trace(\gamma a(x) e_n(x))\nu_g \nonumber \\
&\sim&\sum_{n=0}^{\infty}t^{\frac{n}{2}}\int_{\cls^2}\trace(\gamma a(x) e_n(x))\nu_g. \nonumber 
\end{IEEEeqnarray}
Therefore
\begin{IEEEeqnarray}{rCl}
t^2\trace(\gamma ae^{-t^2D^2})&\sim&\sum_{n=0}^{\infty}t^n\int_{\cls^2}\trace(\gamma a(x) e_n(x))\nu_g \sim 0, \nonumber 
\end{IEEEeqnarray}
as $a(x) e_n(x)$ is a scalar matrix and $\gamma=\left[ {\begin{smallmatrix}
   1  & 0\\
   0 & -1 \\
  \end{smallmatrix}} \right]$. 
  \qed \\
For a pseudo-differential operator $A$, let Wres$(A)$ denote the 
Wodzicki residue  of $A$  (see \cite{Var-2006aa}, page 57).
Let $a^{(1)}=[D,a]$, $a^{(2)}=Fa$ and $a^{(3)}=aF$. 
\bppsn \label{wres1}
For $(j_1,j_2) \neq (1,1)$, one has
\[
\mbox{Wres}(\gamma a_0a_1^{(j_1)}a_2^{(j_2)}\left|D\right|^{-2})= 0.
\]
\eppsn
\prf 
We will show  the result for $j_1=2,j_2=1$. Other cases will follow by similar calculations.
The principal symbol attached to the $\Psi$DO  $D$ and $D^2$  are  $\sigma^D(x,\xi):=\left[ {\begin{smallmatrix}
   0 & \xi \\
   \bar{\xi} & 0 \\
  \end{smallmatrix} } \right]$ and $\sigma^{D^2}(x,\xi):=
	 \left[ {\begin{smallmatrix}
   \left|\xi\right|^2  & 0\\
   0 & \left|\xi\right|^2  \\
  \end{smallmatrix}} \right]$ respectively. 
  Then the principal symbol $\sigma^F(x,\xi)$ of the operator $F=D\left|D\right|^{-1}$ is $\left[ {\begin{smallmatrix}
   0 & \xi \\
   \bar{\xi} & 0 \\
  \end{smallmatrix} } \right]\left[ {\begin{smallmatrix}
   \left|\xi\right|^{-1}  & 0\\
   0 & \left|\xi\right|^{-1}  \\
  \end{smallmatrix} } \right]=\left[ {\begin{smallmatrix}
   0  & \bar{\xi}\left|\xi\right|^{-1}\\
   \xi \left|\xi\right|^{-1} & 0  \\
  \end{smallmatrix} } \right]$. Moreover, since $[D,a]=-ic(da)$ where $c$ denote the clifford action,
  the principal symbol for $[D,a]$ is $-i\left[ \begin{smallmatrix}
   0  & \frac{\partial a}{\partial x_1}-i\frac{\partial a}{\partial x_2}\\
   \frac{\partial a}{\partial x_1}+i\frac{\partial a}{\partial x_2} & 0  
  \end{smallmatrix} \right]$.  Further, as a multiplication operator, $a_0$ is $\Psi$DO of order $0$ with principal
  symbol $\sigma^{a_0}(x,\xi)=\left[ \begin{smallmatrix}
   a_0(x)  & 0\\
   0 & a_0(x)  
  \end{smallmatrix} \right]$. Hence the principal symbol attached to the $\Psi$DO  $\gamma a_0a_1^{(1)}a_2^{(2)}\left|D\right|^{-2}$ is given by
	\[
	\sigma(x,\xi) =
	-i\begin{pmatrix}
   a_0(x)a_1(x)(\frac{\partial a_2}{\partial x_1}+i\frac{\partial a_2}{\partial x_2})\bar{\xi}\left|\xi\right|^{-3}  & 0\\
   0 & -a_0(x)a_1(x)(\frac{\partial a_2}{\partial x_1}-i\frac{\partial a_2}{\partial x_2})\xi\left|\xi\right|^{-3}  
  \end{pmatrix}.
	\] 
For $\left|\xi\right|=1$, $\trace(\sigma(x,\xi))=-2a_0(x)a_1(x)(\xi_1 \frac{\partial a_2}{\partial x_2}-\xi_2 \frac{\partial a_2}{\partial x_2})$.
Hence the Wodzicki residue density of the $\Psi$DO  $\gamma a_0a_1^{(j_1)}a_2^{(j_2)}\left|D\right|^{-2}$ at the point $x$ is
	\[ 
	\int_{\left|\xi\right|=1}\trace(\sigma(x,\xi)) =0.
	\]
Hence
\[
\mbox{Wres}(\gamma a_0a_1^{(j_1)}a_2^{(j_2)}\left|D\right|^{-2})= 0.
\]
This proves the claim.
\qed 
\bppsn \label{wres2} 
The Wodzicki residue of the operator $\gamma a_0da_1da_2\left|D\right|^{-2}$ is given by 
\[
\mbox{Wres}(\gamma a_0da_1da_2\left|D\right|^{-2})= -4\pi i\int_{\cls^2}a_0 da_1 \wedge da_2.
\]
\eppsn
\prf 
 The principal symbol attached to the $\Psi$DO  $\gamma a_0da_1da_2\left|D\right|^{-2}$ is given by
\[
\sigma(x,\xi) = -\begin{pmatrix}
   a_0(x)(\frac{\partial a_1}{\partial x_1}-i\frac{\partial a_1}{\partial x_2})(\frac{\partial a_2}{\partial x_1}
   +i\frac{\partial a_2}{\partial x_2})\left|\xi\right|^{-2}  & 0\\
   0 & -a_0(x)(\frac{\partial a_2}{\partial x_1}
   +i\frac{\partial a_2}{\partial x_2})(\frac{\partial a_2}{\partial x_1}-i\frac{\partial a_2}{\partial x_2})\left|\xi\right|^{-2}  
  \end{pmatrix}.
	\] 
For $\left|\xi\right|=1$, $\trace(\sigma(x,\xi))=-2ia_0(\frac{\partial a_1}{\partial x_1}\frac{\partial a_2}{\partial x_2}-\frac{\partial a_1}
{\partial x_2}\frac{\partial a_2}{\partial x_1})$. Hence
\begin{IEEEeqnarray}{rCl}
\mbox{Wres}(\gamma a_0da_1da_2\left|D\right|^{-2})&=& 
\int_{\cls^2}\left(\int_{\left\|\xi\right\|=1}-2ia_0\left(\frac{\partial a_1}{\partial x_1}\frac{\partial a_2}{\partial x_2}
-\frac{\partial a_1}{\partial x_2}\frac{\partial a_2}{\partial x_1}\right) \sigma \right)dx_1\wedge dx_2 \nonumber \\
&=& -4\pi i\int_{\cls^2}a_0 da_1 \wedge da_2. \nonumber
\end{IEEEeqnarray}
This proves the assertion.
\qed

Let 
$\clb$ be the algebra of pseudo-differential operators of order $0$. Then $\clb$ is a *-subalgebra of $\cll(\clh)$ for which 
all conditions for the WHKAE property of the spectral triple $(C^{\infty}(\cls^2),\clh,D,\gamma)$ holds (see section $4.3$, \cite{ChaSun-2011ab}). 
Let $I_0=I_1=\left\{0\right\}$ and $I_2= \clb$. Then 
$J_0=J_1=\left\{0\right\}\subset J_2=\clb \otimes \cls(\ell^2(\bbn)) \subset J_3=\Sigma^2\clb$ are ideals in $\Sigma^2(C(\cls^2))$ such that
\begin{enumerate}
\item
$a \in J_{\ell} \Longrightarrow d_0(a) \in J_{\ell} \quad \mbox{and} \quad \delta_0(a) \in J_{\ell}$,
\item
for every $a \in J_{\ell}$,  $ t^{\ell}\trace(ae^{-t\left|D_0\right|})$ has an asymptotic expansion near $0$. 
\end{enumerate} 
\blmma \label{phi0}
For $\tilde{a} \in \Sigma^2C^{\infty}(\cls^2)$, one has 
\[
\Sigma^2\phi_0((\gamma \otimes 1)\tilde{a})=0.
\]
\elmma
\prf
 Let $\tilde{a}=a\otimes k$ where $a \in C^{\infty}(\cls^2)$ and $k \in \cls(\ell^2(\bbn))$. Using
  Proposition \ref{sym} and Proposition \ref{compare}, we get
\begin{IEEEeqnarray}{rCl}
t^2\trace((\gamma \otimes 1 )(a\otimes k)e^{-t\left|D_0\right|})=t^2\trace(\gamma ae^{-t\left|D\right|})\trace(ke^{-tN}) \sim 0.
\end{IEEEeqnarray}
It follows from Proposition \ref{res} that  
\begin{displaymath}
 \Sigma^2\phi_0(a\otimes k)= \mbox{Res}_{z=0}(z^{-1}\trace((\gamma a \otimes k)\left|D_0\right|^{-2z}))=0.
\end{displaymath} 
For   $\tilde{a}=1\otimes S^n$ where $n \in \bbz$, we have
\begin{displaymath}
 \trace((\gamma \otimes 1 )(1\otimes S^n)e^{-t\left|D_0\right|})=
\trace(\gamma e^{-t\left|D\right|})\trace(S^ne^{-tN})=0. 
\end{displaymath}
This completes the proof. 
\qed 

\blmma \label{lmmaphi_2}
For $m_1, m_2, m_3 \in \bbn$, one has 
\begin{IEEEeqnarray}{rCl}
\Sigma^2\phi_2(1 \otimes S^{m_1}, 1 \otimes S^{m_2}, 1 \otimes S^{m_3})&=& 0. \nonumber
\end{IEEEeqnarray}
\elmma
\prf
Observe that $d_0(1\otimes S^m)=-mF\otimes S^m$ and $\delta_0^{x}d_0(1\otimes S^m)=-m^{x+1}F\otimes S^m$. Hence we have
\begin{IEEEeqnarray*}{lCl}
\IEEEeqnarraymulticol{3}{l}{\Sigma^2\psi_{x}^{(0)}((\gamma \otimes 1)(1\otimes S^{m_0}),d_0(1\otimes S^{m_1}),d_0(1\otimes S^{m_2})) }\\
&=&  \mbox{Res}_{z=\left(2+\left|x\right| \right)/2}\trace\left((\gamma \otimes S^{m_0})\delta_{0}^{x_1}d_0(1\otimes S^{n_1})
\delta_{0}^{x_2}d_0(1\otimes S^{n_2})\left|D_0\right|^{-2z}\right)\\
&=&  m_1^{x_1+1}m_2^{x_2+1}\mbox{Res}_{z=\left(2+\left|x\right| \right)/2}\trace\left((\gamma \otimes S^{\sum_{i=0}^2m_i})
\left|D_0\right|^{-2z}\right).
\end{IEEEeqnarray*}
Since $\trace\left((\gamma \otimes S^{\sum_{i=0}^2m_i})e^{-t\left|D_0\right|}\right)=
\trace (\gamma e^{-t\left|D\right|})\trace (S^{\sum_{i=0}^2m_i} e^{-tN})= 0$, by Proposition \ref{res}, we have 
$\mbox{Res}_{z=\left(2+\left|x\right| \right)/2}\trace\left((\gamma \otimes S^{\sum_{i=0}^2m_i})
\left|D_0\right|^{-2z}\right)=0$ which further implies that 
$\Sigma^2\psi_{x}^{0}((\gamma \otimes 1)(1\otimes S^{m_0}),d_0(1\otimes S^{m_1}),d_0(1\otimes S^{m_2}))=0$.
Now by applying equation (\ref{e5}),
we get 
\[
\Sigma^2\phi_2(1 \otimes S^{m_1}, 1 \otimes S^{m_2}, 1 \otimes S^{m_3})=0.
\] 
\qed

\bthm \label{thmphi_2}
If for some $i \in \left\{0,1,2\right\}$, $a_i\otimes c_i \in C^{\infty}(\cls^2)\otimes \cls(\ell^2(\bbn))$ then
\begin{IEEEeqnarray*}{rCl}
\Sigma^2\phi_2((\gamma\otimes 1)(a_0\otimes c_0),
                     d_0(a_1\otimes c_1),d_0(a_2\otimes c_2))&=& 
-\frac{i^{\frac{3}{2}}}{\sqrt{2}\pi}\trace(c_0c_1c_2)\int_{\cls^2}a_0 da_1 \wedge da_2.
\end{IEEEeqnarray*} 
\ethm 
\prf
Let $c^{(1)}=c$, $c^{(2)}=Nc$  and $c^{(3)}=-cN$. Then
\[ 
d_0(a\otimes c)= a^{(1)}\otimes c^{(1)}+a^{(2)}\otimes c^{(2)}+a^{(3)}\otimes c^{(3)}.
\] 
Using this and the fact that $a_i\otimes c_i \in C^{\infty}(\cls^2)\otimes \cls(\ell^2(\bbn))$ for some 
$i \in \left\{0,1,2\right\}$, we conclude that the operator 
$(\gamma \otimes 1)(a_0 \otimes c_0)\delta_0^{x_1}d_0(a_1\otimes c_1)\delta_0^{x_2}d_0(a_2\otimes c_2)$ is in the ideal $J_2$. 
Therefore the function 
$ t^2\trace \left((\gamma \otimes 1)(a_0\otimes c_0)\delta_0^{x_1}(a_1\otimes c_1)\delta_0^{x_2}(a_2\otimes c_2)e^{-t\left|D_0\right|}\right) $ has an
asymptotic expansion property near $0$ and hence it follows from equation (\ref{e5}) that  for $x \neq (0,0)$, 
\begin{IEEEeqnarray}{rCl} 
\Sigma^2\psi_x^{(0)}((\gamma \otimes 1)(a_0\otimes c_0),d_0(a_1\otimes c_1),d_0(a_2\otimes c_2))&=&0.\nonumber
\end{IEEEeqnarray}
 For $x=(0,0)$, we have, by Proposition~\ref{wres1},
\begin{IEEEeqnarray*}{lCl}
\IEEEeqnarraymulticol{3}{l}{\Sigma^2\psi_{(0,0)}^{(0)}((\gamma\otimes 1)(a_0\otimes c_0),d_0(a_1\otimes c_1),d_0(a_2\otimes c_2))} \\
&=&\mbox{Res}_{z=1}\trace((\gamma a_0\otimes c_0)d_0(a_1\otimes c_1)d_0(a_3\otimes c_3)\left|D_0\right|^{-2z}) \\
&=&\sum_{j_1=1}^3\sum_{j_2=1}^3 \psi_{(0,0)}^{(0)}(\gamma a_0,a_1^{(j_1)},a_2^{(j_2)})\varphi_0(c_0c_1^{(j_1)}c_2^{(j_2)}) \\
&=& \psi_{(0,0)}^{(0)}(\gamma a_0,da_1,da_2)\trace(c_0c_1c_2).
\end{IEEEeqnarray*} 
 By Connes' Trace Theorem, (see  Theorem $1$, \cite{Con-1988aa}), we have 
\[
\mbox{Res}_{z=1}\trace(A\left|D\right|^{-2z})=2\mbox{Tr}^+(A)=\frac{1}{4\pi^2}\mbox{Wres} A,
\] 
where $A$ is $\Psi$DO  of order $-2$ and  $\mbox{Tr}^+$ denotes the Dixmier trace. 
Note that $\gamma a_0a_1^{(j_1)}a_2^{(j_2)}$ is a $\Psi$DO of order $0$ and hence $\gamma a_0a_1^{(j_1)}a_2^{(j_2)}\left|D\right|^{-2}$ is
$\Psi$DO of order $-2$. 
By Proposition~\ref{wres2}, we get
\begin{IEEEeqnarray*}{lCl}
\IEEEeqnarraymulticol{3}{l}{\Sigma^2\psi_{(0,0)}^{(0)}((\gamma\otimes 1)(a_0\otimes c_0),d_0(a_1\otimes c_1),d_0(a_2\otimes c_2))} \\
&=&2\mbox{Tr}^+(\gamma a_0da_1da_2\left|D\right|^{-2})\trace(c_0c_1c_2) \\
&=& \frac{1}{4\pi^2}\mbox{Wres}(\gamma a_0da_1da_2\left|D\right|^{-2})
                 \trace(c_0c_1c_2)  \\
%& & \qquad \qquad \qquad \qquad \qquad \qquad \qquad \qquad \quad \mbox{by lemma}\quad \ref{wres1}\nonumber \\
&=&  -\frac{i}{\pi}\trace(c_0c_1c_2)\int_{\cls^2}a_0 da_1 \wedge da_2. 
\end{IEEEeqnarray*}
By equation $(\ref{e5})$, we have
\begin{IEEEeqnarray*}{lCl}
\IEEEeqnarraymulticol{3}{l}{\Sigma^2\phi_2((a_0\otimes c_0),d_0(a_1\otimes c_1),d_0(a_2\otimes c_2))} \\
&=&\sum_{x \in \bbn^2}B_{x}^2\Sigma^2\psi_{x}^{(0)}((\gamma\otimes 1)(a_0\otimes c_0),d_0(a_1\otimes c_1),d_0(a_2\otimes c_2))  \\
&=& B_{(0,0)}^2\Sigma^2\psi_{(0,0)}^{(0)}((\gamma\otimes 1)(a_0\otimes c_0),d_0(a_1\otimes c_1),d_0(a_2\otimes c_2))  \\
&=&-\>\frac{i^{\frac{3}{2}}}{\sqrt{2}\pi}\trace(c_0c_1c_2)\int_{\cls^2}a_0 da_1 \wedge da_2. 
\end{IEEEeqnarray*}
This completes the proof. \qed

\subsection{Local index formula for quantum double suspension of non-commutative torus}
Let us recall the definition of non-commutative torus. Throughout we assume that $\theta $ is irrational.
\bdfn
The $C^*$-algebra $A_{\theta}$ is defined as the universal $C^*$-algebra generated by two
unitaries $u$ and $v$ such that $uv=e^{2\pi i\theta}vu$.
\edfn
Define the operators $U$ and $V$ on $\clh:=\ell^2(\bbz^2)$ as follows:
\begin{IEEEeqnarray}{rCl}
Ue_{m,n}&:=&e_{m+1,n} \nonumber \\
Ve_{m,n}&:=&e^{-2 \pi i n\theta}e_{m,n+1} \nonumber 
\end{IEEEeqnarray}
where $\left\{e_{m,n}\right\}$ denotes the standard orthonormal basis of $\ell^2(\bbz^2)$. 
It is well known that $u\mapsto U$ and $v\mapsto V$ gives a faithful representation of the $C^*$-algebra  $A_{\theta}$.

For a function $f(m,n)$ on $\bbz^2$, define the operator $T_f$ as $T_fe_{m,n}:=f(m,n)e_{m,n}$. The group $\bbz^2$ 
acts on the algebra of functions as follows: For $x=(a,b) \in \bbz^2$ and $f(m,n)$, define $x.f:=f(m-a,n-b)$. 
We denote $(1,0)$ by $e_1$ and $(0,1)$ by $e_2$.
Let $\cla_{\theta}$ be the $*$-algebra generated by $u$ and $v$. We consider direct sum representation of
$A_{\theta}$ on $\clh \oplus \clh$. Define $D:= \left[ {\begin{smallmatrix}
   0  & T_{m-in}\\
   T_{m+in} & 0 \\
  \end{smallmatrix} } \right]$ and the grading operator $\gamma:= \left[ {\begin{smallmatrix}
   1  & 0\\
   0 & -1 \\
  \end{smallmatrix} } \right]$.  
	Now note the following commutation relations:
	\begin{IEEEeqnarray}{rCl}
	UT_f&:=&T_{e_1.f}U \nonumber \\
	VT_f&:=&T_{e_2.f}V. \nonumber 
	\end{IEEEeqnarray} 
	It follows that $[T_{m+in},U^{\alpha}V^{\beta}]=(\alpha+i\beta)U^{\alpha}V^{\beta}$ and 
	$[T_{m-in},U^{\alpha}V^{\beta}]=(\alpha-i\beta)U^{\alpha}V^{\beta}$.
	Let us define 
	$\mathcal{D}_k :=\mbox{span}\{T_{P}U^{\alpha}V^{\beta}: \alpha,\beta \in \bbz, P $  is a polynomial of  degree 
	$\leq k \}$ and let $\mathcal{D} := \cup \mathcal{D}_k$. Denote by $\Delta$ the unbounded operator $T_{m^2+n^2}$.
	It follows from proposition $4.13$ \cite{ChaSun-2011ab} and proposition $4.14$ \cite{ChaSun-2011ab} that $(\mathcal{D}, \Delta)$ is a
	differential pair of analytic dimension $2$ with heat kernel expansion property. Now, by proposition $4.15$  \cite{ChaSun-2011ab},
	the spectral triple $(\cla_{\theta},\clh \oplus \clh,D,\gamma)$ has WHKAE property with dimension $2$.
	Consider its quantum double suspension spectral triple 
	$(\Sigma^2(\cla_{\theta}),\clh\otimes \ell^2(\bbn), \Sigma^2D,\gamma\otimes 1)$. 
	It follows from proposition $\ref{p3}$ that its dimension spectrum is contained in $\left\{1,2,3\right\}$ 
	and hence $\Sigma^2\phi_{2n}=0$ for $n>1$. Therefore to compute LIF for this spectral triple, we need to compute 
	$\Sigma^2\phi_{0}$ and $\Sigma^2\phi_{2}$. We will first find asymptotic expansion of $t^2\trace(e^{-t^2\Delta})$.
\begin{IEEEeqnarray}{rCl} 
t^2\trace(e^{-t^2\Delta})&=&2t^2\sum_{m,n \in \bbz}e^{-t^2(m^2+n^2)}. \nonumber \\
&=&2( t\sum_{m\in \bbz}e^{-t^2m^2})(t\sum_{n\in \bbz}e^{-t^2n^2}). \nonumber \\
&=&2\pi \qquad \qquad \mbox{as} \int_{-\infty}^{\infty}{e^{-ax^2}dx}= \sqrt{\frac{\pi}{a}}. \label{normal}
\end{IEEEeqnarray}

\bppsn \label{v} 
For $a_0, a_1,a_2 \in \cla_{\theta}$, one has
 \begin{enumerate}
 \item 
 $\mbox{Res}_{z=1}(\trace(\gamma a_0d(a_1)Fa_2\left|D\right|^{-2z}))= 0.$ 
 \item 
 $\mbox{Res}_{z=1}(\trace(\gamma a_0d(a_1)a_2F\left|D\right|^{-2z}))= 0. $
 \item
 $\mbox{Res}_{z=1}(\trace(\gamma a_0Fa_1d(a_2)\left|D\right|^{-2z}))= 0. $
 \item
 $\mbox{Res}_{z=1}(\trace(\gamma a_0a_1Fd(a_2)\left|D\right|^{-2z}))= 0. $
 \end{enumerate}
\eppsn
\prf
We will show the claim for the first case.  Other cases will follow by similar calculation.
Let $a_i=u^{\alpha_i}v^{\beta_i}$ for $i \in \left\{0,1,2\right\}$.  Then 
\begin{IEEEeqnarray}{rCl}
\left[D,a_1\right]&=& \left( {\begin{array}{cc}
   0 & \left[T_{m-in},U^{\alpha_1}V^{\beta_1}\right]\\
   \left[T_{m+in},U^{\alpha_1}V^{\beta_1}\right] & 0 \\
  \end{array} } \right) \nonumber \\ 
	&=&\left( {\begin{array}{cc}
   0 & (\alpha_1-i\beta_1)U^{\alpha_1}V^{\beta_1}\\
   (\alpha_1+i\beta_1)U^{\alpha_1}V^{\beta_1}& 0 \\
  \end{array} } \right) \nonumber 
	\end{IEEEeqnarray} 
Also, using the identity $\left|D\right|^{-2z}a\left|D\right|^{-2z}\sim \sum_{j=0}^{\infty}{z \choose j} \nabla^j(a)\left|D\right|^{-2j}$, we get

\[
Fa_2=D\left|D\right|^{-1}a_2 \sim \sum_{j=0}^{\infty}{-1/2 \choose j}D \nabla^j(a_2)\left|D\right|^{-1-2j}.
\]
Hence
\[ 
\mbox{Res}_{z=0}\trace(\gamma a_0 da_1 Fa_2 \left|D\right|^{-2-2z})=
\mbox{Res}_{z=0}\trace\left(\sum_{j=0}^{\infty}{-1/2 \choose j}\gamma a_0 da_1 D \nabla^j(a_2)\left|D\right|^{-3-2j-2z}\right).
\]
Since  $ T_{m-in}\Delta^j(a_2)$ and $ T_{m+in}\Delta^j(a_2)$  are elements of
$\mathcal{D}_{j+1}$, $ t^{j+3}\trace\left(\gamma a_0 da_1 D \nabla^j(a_2)e^{-t^2D^2}\right)$ has
asymptotic expansion near $0$. This  implies that 
$\trace(\gamma a_0 da_1 D \nabla^j(a_2)\left|D\right|^{-2z})$ has simple poles 
at $\left\{1/2, \cdots,(j+3)/2\right\}$. Hence for $j \neq 0$,  
$\mbox{Res}_{z=0}\trace(\gamma a_0 da_1 D \nabla^j(a_2)\left|D\right|^{-3-2j-2z})=0.$
 For $j=0$, we have 
	\begin{IEEEeqnarray}{lCl}
	\gamma a_0 da_1 D a_2\left|D\right|^{-3-2z}\nonumber \\
	=\left( {\begin{array}{cc}
  (\alpha_1-i\beta_1)U^{\alpha_0}V^{\beta_0}U^{\alpha_1}V^{\beta_1} T_{m-in}U^{\alpha_2}V^{\beta_2} & 0\\
   0 & - (\alpha_1+i\beta_1)U^{\alpha_0}V^{\beta_0}U^{\alpha_1}V^{\beta_1} T_{m+in}U^{\alpha_2}V^{\beta_2}  \\
  \end{array} } \right)\left|D\right|^{-3-2z} \nonumber \\
	=\left({\begin{array}{cc}
  c_1U^{\alpha_0+\alpha_1+\alpha_2}V^{\beta_0+\beta_1+\beta_2}T_{(m-\alpha_2)-i(n-\beta_2)} & 0\\
   0 & c_2U^{\alpha_0+\alpha_1+\alpha_2}V^{\beta_0+\beta_1+\beta_2}T_{(m-\alpha_2)+i(n-\beta_2)}   \\
  \end{array} }\right)\left|D\right|^{-3-2z} \nonumber 
	\end{IEEEeqnarray}
	where $c_1, c_2$ are constants depending on $\alpha_i, \beta_i$ and $\theta$. 
	It is clear that if $(\alpha_0+\alpha_1+\alpha_2,\beta_0+\beta_1+\beta_2) \neq (0,0)$,
	$\trace(\gamma a_0 da_1 Fa_2\left|D\right|^{-3-2z}) =0$. For $(\alpha_0+\alpha_1+\alpha_2,\beta_0+\beta_1+\beta_2) = (0,0)$, we have
	\begin{IEEEeqnarray}{lCl}
	\mbox{Res}_{z=0}\trace\left(T_{(m-\alpha_2)-i(n-\beta_2)}(T_{\sqrt{m^2+n^2}})^{-3-2z}\right)\nonumber \\
	=\mbox{Res}_{z=0}\sum_{(m,n)\in \bbz^2}\left((m-\alpha_2)-i(n-\beta_2)\right)\left(\sqrt{m^2+n^2}\right)^{-3-2z}\nonumber \\
	=\mbox{Res}_{z=0}\sum_{(m,n)\in \bbz^2}(-\alpha_2+i\beta_2)\left(\sqrt{m^2+n^2}\right)^{-3-2z} \nonumber \\
	=(-\alpha_2+i\beta_2)\mbox{Res}_{z=3/2}\sum_{(m,n)\in \bbz^2}\left(\sqrt{m^2+n^2}\right)^{-2z} \nonumber \\
	=0 . \nonumber 
	\end{IEEEeqnarray}
	Therefore
\begin{IEEEeqnarray}{rCl}
\mbox{Res}_{z=1}\trace(\gamma a_0 da_1 Fa_2 \left|D\right|^{-2z})&=
&\mbox{Res}_{z=0}\trace(\gamma a_0 da_1  Fa_2 \left|D\right|^{-2-2z}). \nonumber \\
&=& \mbox{Res}_{z=0}\trace\left(\sum_{j=0}^{\infty}{-1/2 \choose j}\gamma a_0 da_1 D \nabla^j(a_2)\left|D\right|^{-3-2j-2z}\right). \nonumber \\
&=& 0. \nonumber 
\end{IEEEeqnarray}	\qed

\bppsn \label{x} 
 For $a_0, a_1,a_2 \in \cla_{\theta}$, one has 
 \begin{enumerate}
 \item 
 $\mbox{Res}_{z=1}\trace(\gamma a_0Fa_1Fa_2\left|D\right|^{-2z})=0.$ 
 \item 
 $\mbox{Res}_{z=1}\trace(\gamma a_0Fa_1a_2F\left|D\right|^{-2z})= 0.$
 \item
 $\mbox{Res}_{z=1}\trace(\gamma a_0a_1Fa_2F\left|D\right|^{-2z})=0.$
 \item
 $\mbox{Res}_{z=1}\trace(\gamma a_0a_1F^2a_2\left|D\right|^{-2z})= 0.$
 \end{enumerate}
\eppsn
\prf Let $a_i=u^{\alpha_i}v^{\beta_i}$ for $i \in \left\{0,1,2\right\}$. Note that 
\begin{IEEEeqnarray}{rCl}
Fa_1 &=& D\left|D\right|^{-1}a_1\sim \sum_{j=0}^{\infty}{-1/2 \choose j}D \nabla^j(a_1)\left|D\right|^{-1-2j}. \nonumber 
\end{IEEEeqnarray}
Using this expansion together with the identity
$\left|D\right|^{-2z}a\left|D\right|^{-2z}\sim \sum_{j=0}^{\infty}{z \choose j} \nabla^j(a)\left|D\right|^{-2j}$, we get
\begin{IEEEeqnarray}{lCl}
\mbox{Res}_{z=1}(\trace(\gamma a_0Fa_1Fa_2\left|D\right|^{-2z}))\nonumber \\
=\mbox{Res}_{z=0}(\trace(\gamma a_0Fa_1Fa_2\left|D\right|^{-2-2z})) \nonumber \\
=\mbox{Res}_{z=0}\sum_{j=0}^{\infty}\sum_{l=0}^{\infty}{-1/2 \choose j}{1+j \choose l}\gamma a_0
D\nabla^j(a_1)D\nabla^l(a_2)\left|D\right|^{-4-2j-2l-2z} \nonumber
\end{IEEEeqnarray}
Since $\gamma a_0
D\nabla^j(a_1)D\nabla^l(a_2)$ is element of $\mathcal{D}_{j+l+2}$,
$ t^{j+l+4}\trace\left(\gamma a_0
D\nabla^j(a_1)D\nabla^l(a_2)e^{-t^2D^2}\right)$ has
asymptotic expansion near $0$. This implies that
$\trace(\gamma a_0
D\nabla^j(a_1)D\nabla^l(a_2)\left|D\right|^{-2z})$ has simple poles 
at $\left\{1/2, \cdots,(j+l+4)/2\right\}$. Hence for $(j,l) \neq (0,0)$,  
\[\mbox{Res}_{z=0}\trace(\gamma a_0
D\nabla^j(a_1)D\nabla^l(a_2)\left|D\right|^{-4-2j-2l-2z})=0.\] Therefore
\begin{IEEEeqnarray}{lCl}
\mbox{Res}_{z=1}(\trace(\gamma a_0Fa_1Fa_2\left|D\right|^{-2z}))\nonumber \\
=\mbox{Res}_{z=0}(\trace(\gamma a_0Da_1Da_2\left|D\right|^{-4-2z})) \nonumber\\
=\mbox{Res}_{z=0}\trace \left(\gamma a_0a_1a_2 \left[{\begin{array}{cc}
  T_{m^2+n^2} & 0\\
   0 & T_{m^2+n^2} \\
  \end{array} }\right] \left|D\right|^{-4-2z}\right)\nonumber \\
=e^{-2 \pi i\theta(\alpha_1\beta_0+\alpha_2\beta_0+\alpha_2\beta_1)}\mbox{Res}_{z=0}\trace \left(\left[{\begin{array}{cc}
  U^{\left|\alpha\right|}V^{\left|\beta\right|} & 0\\
   0 & -U^{\left|\alpha\right|}V^{\left|\beta\right|} \\
  \end{array} }\right] \left|D\right|^{-2-2z}\right)\nonumber \\
=0 \nonumber
\end{IEEEeqnarray}
Note that in second step in above calculation, we removed some terms which had zero residues.  Other parts of the claim will follow by this result. 
One can write $\gamma a_0Fa_1a_2F$, $\gamma a_0a_1Fa_2F$ and $\gamma a_0a_1F^2a_2$ 
in the form $\gamma b_0Fb_1Fb_2$. 
and then use this result.  
\qed

\bppsn \label{w}
If $a_i=U^{\alpha_i}V^{\beta_i}$ for $i \in \left\{0,1,2\right\}$ then
\[
\mbox{Res}_{z=1}\trace\left(\gamma a_0da_1da_2\left|D\right|^{-2z}\right)=
\begin{cases}
 0 & \mbox{if}\quad  (\left|\alpha\right|,\left|\beta\right|)\neq 0, \cr 
4\pi i(\alpha_1\beta_2-\alpha_2\beta_1)e^{-2 \pi i\theta(\alpha_1\beta_0+\alpha_2\beta_0+\alpha_2\beta_1)}& \mbox{otherwise} \cr
\end{cases}
\]
\eppsn

\prf 
Note that
\begin{IEEEeqnarray}{lCl}
\gamma da_1 da_2 \nonumber \\ 
=\left[{\begin{array}{cc}
  (\alpha_1-i\beta_1)(\alpha_2+i\beta_2)U^{\alpha_0}V^{\beta_0}U^{\alpha_1}V^{\beta_1}U^{\alpha_2}V^{\beta_2} & 0\\
   0 & -(\alpha_1+i\beta_1)(\alpha_2-i\beta_2)U^{\alpha_0}V^{\beta_0}U^{\alpha_1}V^{\beta_1}U^{\alpha_2}V^{\beta_2} \\
  \end{array} }\right] \nonumber \\
	= e^{-2 \pi i\theta(\alpha_1\beta_0+\alpha_2\beta_0+\alpha_2\beta_1)} \left[{\begin{array}{cc}
  (\alpha_1-i\beta_1)(\alpha_2+i\beta_2)U^{\left|\alpha\right|}V^{\left|\beta\right|} & 0\\
   0 & -(\alpha_1+i\beta_1)(\alpha_2-i\beta_2)U^{\left|\alpha\right|}V^{\left|\beta\right|} \\
  \end{array} }\right] \nonumber 
\end{IEEEeqnarray}
Hence if $(\left|\alpha\right|,\left|\beta\right|)\neq 0$,
$\mbox{Res}_{z=1}\trace\left(\gamma a_0da_1da_2\left|D\right|^{-2z}\right)=0$.  For $(\left|\alpha\right|,\left|\beta\right|)= 0$, we have
\begin{IEEEeqnarray}{rCl}
\mbox{Res}_{z=1}\trace\left(\gamma a_0da_1da_2\left|D\right|^{-2z}\right)&=
& 2i(\alpha_1\beta_2-\alpha_2\beta_1)e^{-2 \pi i\theta(\alpha_1\beta_0+\alpha_2\beta_0+\alpha_2\beta_1)}
\mbox{Res}_{z=1}\trace(\left|\Delta\right|^{-z}), \nonumber \\
&=&4\pi i(\alpha_1\beta_2-\alpha_2\beta_1)e^{-2 \pi i\theta(\alpha_1\beta_0+\alpha_2\beta_0+\alpha_2\beta_1)}.\nonumber \\
 & & \qquad \qquad \qquad \qquad  (\mbox{by equation }  \ref{normal} 
\mbox{ and proposition } \ref{compare}) \nonumber 
\end{IEEEeqnarray} \qed

To compute LIF for quantum double suspension of non-commutative torus,
it is enough to take elements of the form $1\otimes S^n$ and $a \otimes k$ 
where $a=u^{\alpha}v^{\beta}$ and $k \in \cls(\ell^2(\bbn))$. 
\blmma
For $\tilde{a} \in \Sigma^2\cla_{\theta}$, one has
\[
\Sigma^2\phi_0((\gamma \otimes 1)\tilde{a})=0.
\]
\elmma

\prf
Let $\tilde{a}=a\otimes k$ where $a =u^{\alpha}v^{\beta} \in \cla_{\theta}$ and $k \in \cls(\ell^2(\bbn))$. We have
 \[
t^2\trace((\gamma \otimes 1 )(a\otimes k)e^{-t\left|D_0\right|})=t^2\trace(\gamma u^{\alpha}v^{\beta}e^{-t\left|D\right|})\trace(ke^{-tN})\sim 0,
\]
For   $\tilde{a}=1\otimes S^n$ where $n \in \bbz, \trace((\gamma \otimes 1 )(1\otimes S^n)e^{-t\left|D_0\right|})=
\trace(\gamma e^{-t\left|D\right|})\trace(S^ne^{-tN})=0$. It follows from proposition \ref{res} that  $\Sigma^2\phi_0((\gamma \otimes 1)\tilde{a})=
\mbox{Res}_{z=0}(z^{-1}\trace((\gamma \otimes 1)\tilde{a}\left|D_0\right|^{-2z}))=0$. This completes the proof. \qed

\blmma
For $m_1, m_2, m_3 \in \bbn$, one has
\begin{IEEEeqnarray}{rCl}
\Sigma^2\phi_2(1 \otimes S^{m_1}, 1 \otimes S^{m_2}, 1 \otimes S^{m_3})&=& 0. \nonumber
\end{IEEEeqnarray}
\elmma
\prf Proof is exactly same as proof of lemma \ref{lmmaphi_2}.
\qed \\
Let $a^{(1)}=[D,a] ,a^{(2)}=Fa, a^{(3)}=aF$ and $c^{(1)}=c,c^{(2)}=Nc, c^{(3)}=-cN$.
\bthm 
Let $a_i=u^{\alpha_i}v^{\beta_i}$ for $i \in \left\{0,1,2\right\}$.
If for some $i \in \left\{0,1,2\right\}$, $a_i\otimes c_i \in \cla_{\theta} \otimes \cls(\ell^2(\bbn)) $ then 
\begin{IEEEeqnarray}{lCl}
 \Sigma^2\phi_2((\gamma\otimes 1)(a_0\otimes c_0),d_0(a_1\otimes c_1),d_0(a_2\otimes c_2))  \nonumber \\
=
\begin{cases}
 0 & \mbox{if}\quad  (\left|\alpha\right|,\left|\beta\right|)\neq 0, \cr 
2\sqrt{2}\pi i^{\frac{3}{2}}(\alpha_1\beta_2-\alpha_2\beta_1)e^{-2 \pi i\theta(\alpha_1\beta_0+\alpha_2\beta_0+\alpha_2\beta_1)}
\trace(c_0c_1c_2)& \mbox{otherwise} \cr
\end{cases}\nonumber
\end{IEEEeqnarray}

\ethm 
\prf Proof is very similar to proof of theorem \ref{thmphi_2} and we will proceed along the same line. 
By the same reasoning, we get that for $x \neq (0,0)$, 
\begin{IEEEeqnarray}{rCl} 
\Sigma^2\psi_x^{(0)}((\gamma \otimes 1)(a_0\otimes c_0),d_0(a_1\otimes c_1),d_0(a_2\otimes c_2))&=&0.\nonumber
\end{IEEEeqnarray}
 For $x=(0,0)$, we have
\begin{IEEEeqnarray}{lCl}
\Sigma^2\psi_{(0,0)}^{(0)}((\gamma\otimes 1)(a_0\otimes c_0),d_0(a_1\otimes c_1),d_0(a_2\otimes c_2))\nonumber \\
=\mbox{Res}_{z=1}\trace((\gamma a_0\otimes c_0)d_0(a_1\otimes c_1)d_0(a_3\otimes c_3)\left|D_0\right|^{-2z}). \nonumber \\
=\sum_{j_1=1}^3\sum_{j_2=1}^3 \psi_{(0,0)}^{(0)}(\gamma a_0,a_1^{(j_1)},a_2^{(j_2)})\varphi_0(c_0c_1^{(j_1)}c_2^{(j_2)}). \nonumber \\
= \psi_{(0,0)}^{(0)}(\gamma a_0,da_1,da_2)\trace(c_0c_1c_2) \qquad (\mbox{ by proposition } \ref{v}, \ref{x},\ref{w}) \nonumber 
\end{IEEEeqnarray} 
By equation (\ref{e5}), we have
\begin{IEEEeqnarray}{lCl}
\Sigma^2\phi_2((a_0\otimes c_0),d_0(a_1\otimes c_1),d_0(a_2\otimes c_2)) \nonumber \\
=\sum_{x \in \bbn^2}B_{x}^2\Sigma^2\psi_{x}^{(0)}((\gamma\otimes 1)(a_0\otimes c_0),d_0(a_1\otimes c_1),d_0(a_2\otimes c_2)) \nonumber \\
=B_{(0,0)}^2\Sigma^2\psi_{(0,0)}^{(0)}((\gamma\otimes 1)(a_0\otimes c_0),d_0(a_1\otimes c_1),d_0(a_2\otimes c_2)) \nonumber \\
=\begin{cases}
 0 & \mbox{if}\quad  (\left|\alpha\right|,\left|\beta\right|)\neq 0, \cr 
2\sqrt{2}\pi i^{\frac{3}{2}}(\alpha_1\beta_2-\alpha_2\beta_1)e^{-2 \pi i\theta(\alpha_1\beta_0+\alpha_2\beta_0+\alpha_2\beta_1)}
\trace(c_0c_1c_2)& \mbox{otherwise} \cr
\end{cases}\nonumber
 \end{IEEEeqnarray}

This completes the proof. \qed

\section{CONCLUDING REMARKS}
In conclusion we would like to say the following.
\begin{enumerate}
\item So far we have very few instances of actual computations with the LIF. The first demonstration was by Connes (\cite{Con-2004aa}). He considered
the case of quantum $SU(2)$ which is nothing but the quantum double suspension of the circle. Later Pal and Sundar (\cite{PalSun-2010aa}) looked at quantum odd spheres. 
These are nothing but iterated quantum double suspensions of the circle (see example $6.3$, \cite{HonSzy-2002aa}). Here we take up the next higher dimensional case and demonstrate
computation of LIF for the quantum double suspension of two sphere and the noncommutative two torus. According to Connes (\cite{Con-Interview}) 
to understand a concept or result 
one may take a two step procedure. In the first step one makes many computations and then the computations are explained conceptually. Here we are really executing the first step.
Only after a good number of computations we will understand the meaning of individual terms in the formula.

 \item Here in these computations many linear functionals appearing in LIF are zero and evaluation of 
 nonzero functionals involve  only top residues. Next one would probably aim computation of the same for QDS  of $C^*$-algebra of continuous functions on 
 higher dimensional spheres or on surfaces of 
 genus more than one. In these cases,  residues other than the top ones may contribute and will 
 add an extra bit of complications which has to be resolved.
 One approach could be  to compute linear functionals  up to a co-boundary which will remove some unnecessary terms.

 \item
 In case of odd spectral triple, we need the linear functionals $\zeta_{D}^{(m)}$ on $\clb$ and 
 two values $\zeta_{D_0}^{(1)}(1)$ and $\zeta_{D_0}^{(1)}(F\otimes 1)$ or $u_p$ and $v_p$ defined in section $3$ to 
 compute the linear functionals $\zeta_{D_0}^{(m)}$ on $\Sigma^2\clb$. Therefore we can carry this process further and 
 get LIF for iterated QDS of $A$ 
 as soon as we know the two values $\zeta_{D_0}^{(1)}(1)$ and $\zeta_{D_0}^{(1)}(F\otimes 1)$ at each stage of iteration.
 Similarly one can iterate the process and 
 get LIF in case of even spectral triple also.
 \item
 We used WHKAE property very crucially in our calculations. Can we drop this property and still compute LIF
 for QDS spectral triple ?
 
 \item In another direction one could pursue investigations along the lines of (\cite{ConMos-2014},\cite{FK}). So far noncommutative tori at various dimensions were the only tractable examples where one could compute with the machinery of local index formula. Now with these models are also available for similar analysis and one would like to investigate them along these lines.

 \end{enumerate}

% \noindent\begin{footnotesize}\textbf{Acknowledgement}:
% The second author would like to thank Prof. Arup Kumar Pal, his supervisor, for his support. He would also like to thank S. Sundar and Satyajit Guin
% for their help during his visit to Institute of Mathematical Sciences, Chennai.
% \end{footnotesize}

%	
	%\bibliography{references}
%\bibliographystyle{plain}
%	

\end{document}